\documentclass{article}
\usepackage{amssymb}
\usepackage{amsfonts}
\usepackage{euscript}

\begin{document}

\vskip2cm
\newcommand{\gothS}{\digamma}
\newcommand{\Varm}{\EuScript{M}}
\newcommand{\D}{\EuScript{D}}
\newcommand{\B}{\EuScript{B}}
\newcommand{\T}{\EuScript{T}}
\newcommand{\gothD}{\mho}
\newcommand{\rl}{{\bf R}^1}
\newcommand{\rr}{{\bf R}}
\newcommand{\rk}{{\bf R}^k}
\newcommand{\rp}{{\bf R}^p}
\newcommand{\rn}{{\bf R}^n}
\newcommand{\ri}{{\bf R}^i}
\newcommand{\nn}{{\bf N}}
\newcommand{\zplus}{{\bf Z}_{+}}
\newcommand{\zz}{{\bf Z}}
\newcommand{\rntn}{{\bf R}^{n \times n}}
\newcommand{\rnextnex}{{\bf R}^{n \times n}}
\newcommand{\rnextpex}{{\bf R}^{n \times p}}
\newcommand{\rnex}{{\bf R}^{n }}
\newcommand{\rpex}{{\bf R}^{p}}
\newcommand{\rktk}{{\bf R}^{k \times k}}
\newcommand{\Linfty}{L_{\infty}(I;\rl)}
\newcommand{\Ll}{L_1(I;\rl)}
\newcommand{\Cheb}{C(I;\rl)}
\newcommand{\Chebk}{C(I;\rk)}
\newcommand{\Gchebk}{C^1(I;\rk)}
\newcommand{\Gcheb}{C^1(I;\rl)}
\newcommand{\zt}{\zeta}
\newcommand{\lmb}{({\lambda}_1,...,{\lambda}_k{)}^T}
\newcommand{\nrho}{\rho}
\newcommand{\neta}{\eta}
\newcommand{\rbn}{{\bf R}^{N}}
\newcommand{\xov}{{\tilde x}^{0}}
\newcommand{\tauex}{{\overline \tau }}
\newcommand{\ff}{{F }}
\newcommand{\rtwo}{{\bf R}^2}
\newcommand{\rthree}{{\bf R}^3}
\newcommand{\rfour}{{\bf R}^4}
\newcommand{\rmnu}{{\bf R}^{m_{\nu}}}
\newcommand{\rntm}{{\bf R}^{n \times m}}
\newcommand{\rmiplusl}{{\bf R}^{m_{i+1}}}
\newcommand{\rmlplusmi}{{\bf R}^{m_1+ \ldots +m_i}}
\newcommand{\rmnuplusl}{{\bf R}^{m_{\nu + 1}}}
\newcommand{\rkplusl}{{\bf R}^{k + 1}}
\newcommand{\rmlplusdp}{{\bf R}^{m_1+ \ldots +m_p}}
\newcommand{\rmpplust}{{\bf R}^{m_{p+2}}}

\newcommand{\rml}{{\bf R}^{m_1}}
\newcommand{\rmt}{{\bf R}^{m_2}}
\newcommand{\rmpmin}{{\bf R}^{m_{p-1}}}
\newcommand{\rsumpmin}{{\bf R}^{m_1 + ...+m_{p-1}}}
\newcommand{\rmtm}{{\bf R}^{m \times m}}
\newcommand{\rmtrmpl}{{\bf R}^{m \times m_{p+1}}}
\newcommand{\rmp}{{\bf R}^{m_p}}
\newcommand{\rmm}{{\bf R}^{m}}
\newcommand{\rmpplusl}{ {\bf R}^{m_{p+1}} }
\newcommand{\rmppluslzpt}{ {\bf R}^{m_{p+1}} }
\newcommand{\rmpplusldt}{ {\bf R}^{m_{p+1}} }
\newcommand{\rmitmj}{ {\bf R}^{m_{i} \times m_{j}} }
\newcommand{\rktrmpl}{{\bf R}^{k \times m_{p+1}}}
\newcommand{\rktkplusl}{{\bf R}^{k \times (k+1)}}
\newcommand{\rmltmtml}{ {\bf R}^{m_{1} \times {(m_{2}-m_1)} }   }
\newcommand{\rmltml}{ {\bf R}^{m_{1} \times m_{1} }   }
\newcommand{\rmtminml}{ {\bf R}^{m_{2}-m_{1}} }
\newcommand {\rmminmp}{ {\bf R}^{m-m_{p}} }
\newcommand {\rrmmpp}{ {\bf R}^{ ({m-m_{p}})\times m_{p} }}
\newcommand{\rmtminm}{{\bf R}^{m_{2} -m} }
\newcommand {\rkminmp}{ {\bf R}^{k-m_{p}} }
\newcommand {\rrkmpp}{ {\bf R}^{ ({k-m_{p}})\times m_{p} }}
\newcommand{\rmtmink}{{\bf R}^{m_{2} -k} }
\newcommand{\rmi}{{\bf R}^{m_i}}
\newcommand{\rmj}{{\bf R}^{m_j}}
\newcommand{\rktm}{{\bf R}^{k \times m}}
\newcommand{\nxi}{\xi}
\newcommand{\rmlplusmp}{{\bf R}^{m_1 + ...+m_{p}}}

\newcommand{\ekk}{E^{2k+2}}

\newcommand{\mll}{\parallel   w   (\cdot)   {\parallel}_{C ([t_0,T];\rpex) }}
\newcommand{\ml}{{\parallel   w   (\cdot)   {\parallel}_{C ([t_0,T];\rpex) }^{-1}}}
\newcommand{\dlm}{{\delta \lambda} = ( {\delta \lambda}_1,\ldots,{\delta\lambda}_k) \in { \bf R}^k}
\newcommand{\wdlms}{w(s)}
\newcommand{\wdlmt}{w(\tau)}

\newcommand{\uhlms}{{\hat u}(s)+{\sum_{j=1}^{k}{ {\lambda_j}{\hat w}_j(s)}}}
\newcommand{\uhlmt}{{\hat u}(\tau)+{\sum_{j=1}^{k}{ {\lambda_j}{\hat w}_j(\tau)}}}

\newcommand{\uhtlms}{{\hat u}_{\lambda}(\cdot)}
\newcommand{\uhtlmt}{{\hat u}_{\lambda}(\tau)}

\newcommand{\wblms}{{\sum_{j=1}^{k}{ {\bar {\lambda_j}}{w}_j(s)}}}
\newcommand{\wblmt}{{\sum_{j=1}^{k}{ {\bar {\lambda_j}}{w}_j(\tau)}}}

\newcommand{\whblms}{{\sum_{j=1}^{k}{ {\bar {\lambda_j}}{\hat {w}_j}(s)}}}
\newcommand{\whblmt}{{\sum_{j=1}^{k}{ {\bar {\lambda_j}}{\hat {w}_j}(\tau)}}}

\newcommand{\wmbblmt}{\hat {w_{\bar{\lambda}}}(\cdot)}
\newcommand{\whbblms}{\hat {w_{\bar{\lambda}}}(s)}
\newcommand{\whbblmt}{\hat {w_{\bar{\lambda}}}(\tau)}
\newcommand{\whbblm}{\hat {w_{\bar{\lambda}}}}

\newcommand{\dfdx}{\frac {\partial f}{\partial x}}
\newcommand{\dfidx}{\frac {\partial {f_i}}{\partial x}}
\newcommand{\dgdx}{\frac {\partial g}{\partial x}}
\newcommand{\dgidx}{\frac {\partial {g_i}}{\partial x}}

\newcommand{\dfdu}{\frac {\partial f}{\partial u}}
\newcommand{\dfidu}{\frac {\partial {f_i}}{\partial u}}
\newcommand{\dgdu}{\frac {\partial g}{\partial u}}
\newcommand{\dgidu}{\frac {\partial {g_i}}{\partial u}}

\newcommand{\ulms}{\overline v (s) }
\newcommand{\ulmt}{\overline v (\tau)  }
\newcommand{\udlms}{{\overline v} (s) + w (s)}
\newcommand{\udlmt}{\overline v (\tau ) + w (\tau)}
\newcommand{\utdlms}{\overline v (s) + {\bar {\theta}_i} (s)  w(s)}
\newcommand{\utdlmt}{\overline v (\tau) + {\bar   {\vartheta}_i }  (s,\tau) w(\tau)}

\newcommand{\xlms}{x_{\overline v (\cdot )}(x^0,  s)}
\newcommand{\xlmt}{x_{\overline v (\cdot )}(x^0,   {\tau})}

\newcommand{\xhlm}{\hat {x_{\lambda}}}

\newcommand{\xhblms}{\hat {x_{\bar{\lambda}}}(s)}
\newcommand{\xhblmt}{\hat {x_{\bar{\lambda}}}({\tau})}
\newcommand{\xhblm}{\hat {x_{\bar{\lambda}}}}
\newcommand{\xdlms}{x_{\overline v (\cdot ) + w(\cdot)}( x^0 ,s)}
\newcommand{\xdlmt}{x_{\overline v (\cdot ) + w(\cdot)}( x^0 ,\tau)}
\newcommand{\xlm}{x_{\overline v (\cdot)} }
\newcommand{\xdlm}{x_{ \overline v (\cdot) + w(\cdot) } }
\newcommand{\xilm}{x_{\overline v (\cdot) }^i}
\newcommand{\xidlm}{x_{\overline v (\cdot)  + w(\cdot) }^i}
\newcommand{\xtdlms}{x_{\overline v (\cdot) } (x^0,s) + {{\theta}_i(s)} (x_{\overline v (\cdot) + w(\cdot) }(x^0,s) - x_{\overline v (\cdot)}(x^0,s))}
\newcommand{\xtdlmt}{x_{\overline v (\cdot)}(x^0,\tau) + {{\vartheta}_i(s,\tau )}(x_{\overline v (\cdot) + w (\cdot)}(x^0,\tau) - x_{\overline v (\cdot)}(x^0,\tau))}
\newcommand{\zidlmlm}{z_{{w(\cdot)},{\overline v (\cdot)}}^i}
\newcommand{\zdlmlm}{z_{{w(\cdot)},{\overline v (\cdot)}}}
\newcommand{\zbdlmlm}{z_{{\bar {\lambda}},{\lambda +{\delta \lambda}}}}
\newcommand{\zblm}{z_{{\bar {\lambda}},{\lambda}}}


\begin{large}

\newcommand{\n}[1]{\refstepcounter{equation}\label{#1}
\eqno{(\arabic{equation})}}
\renewcommand{\theequation}{\arabic{equation}}

\begin{Large}
\begin{center}
{\bf A generalized triangular form and its global controllability}
\end{center}

$$ \; $$

\begin{center}
 Svyatoslav S. Pavlichkov
\end{center}

\end{Large}
$$\; $$

$$ \; $$
\begin{large}
\begin{center}
{\bf Abstract.}
\end{center}
We investigate a new class of nonlinear control systems of O.D.E.,
which are not feedback linearizable in general. Our class is a
generalization of the well-known feedback linearizable systems,
and moreover it is a generalization of  the triangular (or
pure-feedback) forms studied before. The definition of our class
is global, and coordinate-free, which is why the problem of the
equivalence is solved for our class in the whole state space at
the very beginning. The goal of this paper is to prove the global
controllability of our nonlinear systems. We propose to treat our
class as a new canonical form which is a nonlinear global analog
of the Brunovsky canonical form on the one hand, and is a global
and coordinate-free generalization of the triangular form on the
other hand.

\end{large}
$\qquad $

{\bf Key words:} Nonlinear control, triangular form, global
controllability, feedback linearization.

{\bf Mathematics subject classification:} 93C10, 93B10,
93B11, 93B05, 93B52.

$\qquad$
\begin{center}
{\bf 1. Introduction. }
\end{center}

One of the most important problems in the nonlinear control theory
is the problem of classification of nonlinear control systems of
O.D.E., i.e., that of finding the transformation of a nonlinear
control system into its simplest canonical form along with finding
such canonical forms by using some invariants which do not depend
on the choice of local coordinates. Beginning with
\cite{jakubczyk},\cite{hunt}, a complete theory of feedback
linearization was created --
\cite{cheng_isidori},\cite{krishchenko1}, \cite{krishchenko5},
\cite{gardner1}, \cite{gardner2},\cite{D_Andrea}, \cite{MRicc},
etc. However, even some simple mechanical systems do not satisfy
the Respondek-Jakubczyk-Hunt-Su-Meyer conditions; in addition, the
concept of feedback linearization is essentially local. This
inspired many authors to further investigations and to attempts to
generalize the feedback linearization theory.

One possible approach, which is is very popular and has various
applications, is the concept of differential flatness
\cite{fliess}. However, this notion is as local as that of
feedback linearization, and moreover no general criterion of
differential flatness has been obtained.

Another way is to deal with the triangular, or pure-feedback form
instead of the Brunovsky canonical form. Triangular systems were
introduced in \cite{kor1} as early as 1973 (i.e., even before
\cite{jakubczyk},\cite{hunt}) as a first example of a nonlinear
system which is feedback linearizable. Nevertheless, a triangular
system
$$ \left\{ \begin{array}{l}
    \dot z_i = f_i (z_1,...,z_{i{+}1}), \;\; \; i=1,...,n-1;  \\
   \dot z_n = f_{n} (z_1,...,z_n,v);
\end{array}\right.  $$
 is feedback linearizable only in the so-called "regular" case,
 i.e., when the conditions of regularity $|\frac{\partial f_i}{\partial x_{i{+}1}}| {\not=} 0, $
 $i{=}1,...,n$ hold; otherwise (the "singular case"), the system
 does not satisfy the Respondek-Jakubczyk-Hunt-Su-Meyer conditions, and,
 therefore, is not feedback linearizable. The singular case was
 investigated by Respondek in 1986 (see \cite{Resp}), and by Celikovsky and
 Nijmeijer in 1996 (see \cite{Cel_Nij}). In these works, the triangular systems
 are studied under the assumption that the set of the regular
 points is open and dense in the whole state space, however. This is not true
 even for some simple examples (see, for instance system (\ref{ex9}) from
 the current paper).

 That is why, we want to find and to investigate a
 generalization, of the triangular form, which contains all the
 previous triangular forms studied before (including the singular
 case) on the one hand, and which can be investigated globally
 (including the problems of controllability, stabilization,
 feedback equivalence, etc.) on the other hand. We offered such a
 generalization in \cite{ssp_op2}, and solved completely the problem of global
 robust controllability for this class (moreover, the controls constructed
 were actually closed-loop to some extent). However, the problem
 of global equivalence of a control system to a system from \cite{ssp_op2}
 remained open. In this work, we introduce a generalization of the
 triangular systems considered in \cite{ssp_op2}, in global coordinate-free
 terms. The main goal of the current paper is to prove that our
 generalized triangular form is globally controllable.

In the future, we want to investigate in more detail the
relationship between the triangular form from \cite{ssp_op2} and
the class from the current paper. As we can see from example 3.1,
the class of "generalized triangular form" is wider than that from
\cite{ssp_op2}. On the other hand, the construction of example 3.1
is based on triangular system (\ref{ex9}). To what extent our
generalized triangular form can be reduced to the triangular form
in the singular case remains an open question.

{\bf Acknowledgements.} A part of this work was complete when the
author was visiting the Department of Mathematics, Louisiana State
University. The author is grateful to Professor Jimmie D. Lawson
for having organized the visit, for his hospitality, and for
discussion about these results.

\vskip5mm

$\qquad$
\begin{center}
{\bf 2. Notation, and preliminaries. }
\end{center}

Let $\Varm$ be a smooth manifold of dimension $n$, and $x\mapsto
v(x)$ be a smooth vector field on $\Varm.$ In general, $v(\cdot)$
can be defined on some open subset of $\Varm$ only; next we denote
this subset by $\D_{v}.$ Let $x^0$  be in $\Varm.$ By $I{\ni}
t{\mapsto}\Phi_{v}^{t}(x^0) $ we denote the (maximal) integral
curve $t\mapsto x(t)$ of $\dot x  = v(x)$ with $x(0) = x^0.$ Of
course, for each $t \in I$ the map $x \mapsto \Phi_v^t(x)$ is (at
least) a diffeomorphism of some neighborhood of $x^0$ onto some
neighborhood of $\Phi_v^{t} (x^0)$ (and, if, for some $s \in I,$
$\Phi_v^s (x)$ is well-defined for all $x \in \D_{v},$ then $x
\mapsto \Phi_v^s(x)$ is a global diffeomorphism of $\D_v$ onto
$\D_v$)

For every fixed $t\in I,$ every $x$ in a neighborhood of $x^0$,
and every $\xi \in T \Varm_x $ by $\left( \Phi_{v}^{t}
\right)_{\ast} \xi$ we denote the image of $\xi$ under the tangent
map of the diffeomorphism $y\mapsto \Phi_{v}^{t} (y)$ at point
$x.$ (Actually, $\left( \Phi_{v}^{t} \right)_{\ast} \xi$ depends
on two arguments $\xi$ and $x,$ and we should write $\left(
\Phi_{v}^{t} \right)_{\ast} (x,\xi),$ in general, but in our case
it will be always clear at which point $x\in\Varm$ we consider the
tangent map, which is why we write $(\Phi_{v}^{t})_{\ast}\xi$
without any ambiguity.)

In addition, if $V$ is a vector space, then, for $A\subset V,$ and
$B\subset V,$ we denote by $A+B$ the set $\{ x+y \; | \; \; x\in
A, \; y\in B \}$ (in our situation $V$ will be $T{\Varm}_x$ for
some smooth manifold $\Varm$ and some $x \in \Varm$).

If $\Delta(\cdot)$ is a smooth integrable distribution on $\Varm,$
(which means that the dimension $\dim \Delta(x)$ equals $k \le n$
for some fixed $k=1,...,n,$ and for all $x \in \Varm,$ and $\Delta
(\cdot)$ is involutive at each point $x \in \Varm$) then, for each
$x^0 \in \Varm,$ we can consider its orbit, or the maximal
integral manifold $\Varm_{\Delta} (x^0)$ defined as the set of all
points $y \in \Varm$ given by
$$ y = \left( \Phi_{v_1}^{t_1} \circ \Phi_{v_2}^{t_2} \circ ... \circ \Phi_{v_N}^{t_N} \right) (x^0) \n{r1}$$
with arbitrary $N\ge 1,$ arbitrary $t_i \in \rr,$ $i=1,...,N,$ and
arbitrary smooth vector fields $v_i(\cdot)$ such that, for every
$i=1,...,N,$ and every $x \in \D_{v_i}$ we have $v_i(x) \in
\Delta(x).$ Also we will use a more brief form of (\ref{r1}):
$$y=\Phi_{v}^{T} (x^0) \; \; \; \; \; \; \; \mbox{ with }
T=(t_1,...,t_N), \; \; v=(v_1,...,v_N).\n{r1a}$$ By $\Phi_{v}^{-T}
(\cdot)$ we denote the inverse diffeomorphism, i.e. $ \left(
\Phi_{v_N}^{-t_N} \circ \Phi_{v_{N-1}}^{-t_{N-1}} \circ ... \circ
\Phi_{v_1}^{-t_1} \right) (\cdot)$

We write by definition $v(\cdot) \in \Delta (\cdot)$ iff $v(x) \in
\Delta(x)$ for each $x\in \D_v.$ Let us recall that, if
$v(\cdot),$ and $w(\cdot)$ are smooth vector fields defined on
some open subset $\D \subset \Varm,$ then by $[v,w](\cdot)$ we
denote their Lie bracket defined (in any coordinates) as $[v,w](x)
= \frac{\partial w}{\partial x} v -\frac{\partial v}{\partial
x}w.$ Finally, for $A\subset \Varm,$ we denote by $\overline{A}$
the closure of $A$ in $\Varm.$

\begin{center}
{\bf 3. Main result. }
\end{center}

We consider a control system
$$\dot x = a(x)+ \beta(x,u)b(x) \n{r3}$$
with states $x\in \Varm, $ and controls $u\in \rl,$ where $\Varm$
is a simply connected manifold, $a(\cdot),$ $b(\cdot),$ are smooth
vector fields (of class $C^{n{+}1}$ at least) on $\Varm,$ and
$\beta(\cdot,\cdot)$ is a smooth (of class $C^{n{+}1}$) scalar
function on $\Varm.$ Next we suppose that $\Varm=\rn$ just to make
the arguments clearer, however our technique works for arbitrary
simply connected manifold $\Varm.$ We assume that $a(\cdot),$
$b(\cdot),$ and $\beta (\cdot,\cdot)$ satisfy the following
conditions

{\bf (A)} {\it For each $x\in \Varm,$ we have $b(x) \not= 0,$ and
$\beta(x,\rl)=\rl.$ In other words, the set
$\Delta_0(x):=\{\beta(x,u)b(x)\;|\; \; u \in \rl \}$ is a
1-dimensional subspace of  $T \Varm_x {=}\rn$ for each $x\in
\Varm.$}

(Of course, the distribution $x\mapsto \Delta_0(x)$ is integrable
in the whole $\Varm, $ and, for each $x^0 {\in} \Varm,$ the
corresponding maximal integral manifold of $\Delta_0(\cdot)$ is
 the (maximal) trajectory $t\mapsto \Phi_b^t(x^0)$).

{\bf (B) } {\it  Let $k$ be in $\{1,...,n{-}1 \}.$ Assume that the
distribution $\Varm{\ni}x{\mapsto}\Delta_{k-1}(x) \subset T
\Varm_x {=}{\rn}$ is already  constructed: for $k{=}1,$ see
condition (A), for $k{\ge}1,$  we define $\Delta_k (\cdot)$ by
induction as below. Then:}

{\bf (B1)} {\it We require that $x \mapsto  \Delta_{k{-}1} (x)$ is
of rank $k$ for every $x{\in}\Varm,$ and is involutive at each
$x{\in}{\Varm}.$}

{\bf (B2)} {\it Given $x^0{\in}{\Varm},$ by  $\Delta_k(x^0)$
denote the set
$$ \Delta_k(x^0) := \Delta_{k{-}1} (x^0) + \{ (\Phi_v^T)_{\ast} a(\Phi_v^{-T} (x^0)) - a(x^0) \; | \; \; \forall
N{\ge}1 \; \; \;  \forall T{=}(t_1,...,t_N)  $$ $$ \forall
v(\cdot){=}(v_1(\cdot),...,v_N(\cdot)) \; \;  {\rm such \; that}
\; v_i(x){\in}\Delta_{k{-}1}(x) \;  {\rm for \; all} \;
 x{\in}{\D}_{v_i}, \; i{=}1,...,N   \} \n{r4}$$
 We require that, for each
$x^0{\in}\Varm,$ the set $\Delta_k(x^0)$ is a $(k{+}1)$ -
dimensional subspace of $T \Varm_{x^0} = \rn$ and that the
distribution $x \mapsto \Delta_k (x)$ is involutive (for all
$k=0,\ldots,n-1,$ as we mentioned before.)}

We emphasize that, for each fixed $x^0,$ we obtain $\Delta_k(x^0)$
in (\ref{r4}) by taking all admissible $v_i(\cdot)$ from
$\Delta_{k{-}1}$ (i.e., along the maximal integral manifold
$\Varm_{\Delta_{k{-}1}}(x^0)$ defined in (\ref{r1})).

Let us remark that conditions (A), (B1), (B2) are global analog of
the conditions from  \cite{jakubczyk}, \cite{hunt}.

If a smooth system $\dot x = f(x,u)$ is locally feedback
equivalent to the triangular form, then (see \cite{kor2})
$f(\cdot,\cdot)$ have (locally) the form (\ref{r3}):
$f(x,u)=a(x)+\beta(x,u)b(x)$ with some smooth vector fields
$a(\cdot),$ $b(\cdot),$ $b(x)\not=0$ and with  some smooth scalar
function $\beta(\cdot,\cdot),$ regardless of whether this
triangular form satisfies the regularity conditions
$|\frac{\partial f_i}{\partial x_{i{+}1}}| \not=0,$ or we deal
with the singular case.

Furthermore, any triangular system
$$ \left\{ \begin{array}{l}
    \dot x_i = f_i (x_1,...,x_{i{+}1}), \;\; \; i=1,...,n-1;  \\
   \dot x_n = f_{n} (x_1,...,x_n,u);
 \end{array}\right. \; \; \;\; \; x=(x_1,...,x_n{)}^T \in \rn, \; \; u\in\rl \n{r5} $$
such that $f_i$ are smooth, and $f_i(x_1,...,x_i,\rl) =\rl,$ for
all $i=1,...,n$  and all $(x_1,...,x_i) \in \ri$ (see
\cite{ssp_op2}) satisfies our conditions (A),(B1),(B2)

Conversely, assume that system (\ref{r3}) satisfies (A),(B1),(B2).
Pick any $x^0 \in \Varm,$ and let  $\zeta = \varphi(x) : U(x^0)
{\subset} \Varm \rightarrow V(\varphi(x^0)) {\subset} \rn$
($\zeta_i=\varphi_i(x),$ $i=1,...,n$) be a diffeomorphism of a
neighborhood $U(x^0)$ of $x^0$ onto a neighborhood
$V(\varphi(x^0))$ of $\varphi (x^0).$

{\bf Definition 3.1} {\it We say that coordinates $\zeta_i$
 are {\bf canonical} for system (\ref{r3}), or the map
$x{\mapsto}\zeta{=}\varphi(x)$ defines canonical coordinates for
system (\ref{r3}) (or canonical coordinates for the corresponding
sequence of nested regular integrable distributions
$\Delta_0(\cdot),...,\Delta_{n{-}1} (\cdot)$) in $U(x^0),$ iff,
for each $k=0,...,n{-}2,$ the set $\varphi^{-1} (D_k)$ with $$D_k
:= \{ (\zeta_1,...,\zeta_n) \in V(\varphi(x^0)) \; | \; \; \zeta_i
= {\rm const}, \; \; i = 1,...,n{-}k{-}1 \}$$ is an integral
manifold of $\Delta_k$ in $U (x^0).$}

Equivalently, coordinates $y_i{=}\varphi_i(x)$ $i=1,...,n,$ are
canonical for system (\ref{r3}) in $U(x^0)$ with small enough
$U(x^0),$ iff $$\Delta_k (x) = \{ \xi {\in} T \Varm_x {=} \rn \; |
\; \; <\nabla \varphi_i (x), \xi> = 0, \; \; i = 1,...,n{-}k{-}1
\}, \; \; \; \; k=0,1,\ldots,n-2.$$

{\bf Remark 3.1.} It is easy to prove that, if coordinates
$\zeta_i{=}\varphi_i(x)$ are canonical for (\ref{r3}) in a
neighborhood of some $x^0{\in}{\Varm},$ then, ({\bf locally!}) in
some neighborhood of $x^0,$ this change of coordinates
$\zeta_i{=}\varphi_i(x)$ brings the dynamics of (\ref{r3}) to the
following triangular from
$$ \left\{
\left(\begin{array}{l}
    \dot \zeta_1\\
    \ldots\\
    \dot \zeta_{n{-}2}\\
    \dot \zeta_{n{-}1}\\
    \dot \zeta_n \end{array} \right)
= \left(\begin{array}{c}
    f_1 (\zeta_1,\zeta_2)\\
    \ldots\\
    f_{n{-}2} (\zeta_1,...,\zeta_{n{-}1})\\
    f_{n-1} (\zeta_1,...,\zeta_{n})\\
    f_{n} (\zeta_1,...,\zeta_n) \end{array} \right) +
    \tilde{\beta} (\zeta,u) \left(\begin{array}{c}
    0\\
    \ldots\\
    0\\
    0\\
    g_{n} (\zeta_1,...,\zeta_n) \end{array} \right)
 \right.  \n{rast} $$
where $g_n(\zeta){\not=}0$  in some neighborhood of $\varphi(x^0)$
(but, of course the regularity conditions $|\frac{\partial
f_i}{\partial \zeta_{i{+}1}}| \not= 0$ do not hold, in general).
However, this is true only locally, as we can learn from the
following example.

{\bf Example 3.1} Consider the vector fields  $v_1(\cdot)$ and
$v_2(\cdot)$ in $\rtwo$ given by $v_1(x) = (\cos x_1, - \sin
x_1{)}^T, $ and $v_2 (x) = (\sin x_1, \cos x_1{)}^T.$  Put $x_N^0
= (\pi N, 0) \in \rtwo$ for all $N \in {\bf Z}.$ For each $N \in
{\bf Z},$  the map given by $\rtwo \ni (z_1,z_2) \mapsto
(\Phi_{v_2}^{z_2} \circ \Phi_{v_1}^{z_1}) (x_N^0)$ is a
diffeomorphism of $\rtwo$ onto $]\pi (N-1), \pi(N+1) [ \times \rr$
(we do not solve the corresponding differetial equations
explicitely just because the solution can not be written out as a
combination of standard functions). Fix any $z_1>0.$ Denote by
$A(z_1) \in \rtwo$ the intersection of the trajectory $z_2\mapsto
(\Phi_{v_2}^{z_2} \circ \Phi_{v_1}^{z_1}) (x_0^0)$ with $\{
\frac{\pi}{2} \} \times \rr ;$  then from the symmetry of the
curve $z_2 \mapsto (\Phi_{v_2}^{z_2} \circ \Phi_{v_1}^{z_1})
(x_0^0) $ w.r.t. $\{ \frac{\pi}{2} \} \times \rr$ we obtain that
there is a unique solution $({\tilde z}_2 (z_1), {\tilde z}_1
(z_1))$ of $(\Phi_{v_2}^{{\tilde z}_2 (z_1)} \circ
\Phi_{v_1}^{z_1}) (x_0^0) = \Phi_{v_1}^{{\tilde z}_1 (z_1)}
(x_1^0)$ w.r.t. $({\tilde z}_2 (z_1),{\tilde z}_1 (z_1))$;
$$ {\tilde z_1} (z_1) = z_1; \; \; \; \; \; \; \; \; A(z_1) = \left(\Phi_{v_2}^{\frac{{\tilde z}_2 (z_1)}{2}} \circ \Phi_{v_1}^{z_1}\right) (x_0^0); \n{ex6}$$
and there is a unique solution $(z_2^{\ast} (z_1), z_1^{\ast}
(z_1))$ of the nonlinear equation $\left( \Phi_{v_2}^{z_2^{\ast}
(z_1)} \circ \Phi_{v_1}^{z_1} \right) (x_1^0) =
\Phi_{v_1}^{z_1^{\ast} (z_1)} (x_0^0)$ w.r.t $(z_2^{\ast} (z_1),
z_1^{\ast} (z_1)),$ and $$z_1^{\ast} (z_1) = z_1; \; \; \; \;
z_2^{\ast} (z_1) = - {\tilde z}_1(z_1);   \; \; \; \; \; A(z_1) =
 \left(\Phi_{v_2}^{\frac{{z}_2^{\ast} (z_1)}{2}} \circ
\Phi_{v_1}^{z_1}\right) (x_1^0); \n{ex7}$$ (in addition, the image
of the curve $\rr\ni z_2 \mapsto (\Phi_{v_2}^{z_2}
\circ\Phi_{v_1}^{z_1}) (x_1^0)$ coincides with the image of
$\rr\ni z_2 \mapsto (\Phi_{v_2}^{z_2} \circ\Phi_{v_1}^{z_1})
(x_0^0)$).

Similarly, for any fixed $z_1<0,$ there is a unique solution
$({\hat z_2} (z_1), {\hat z}_1 (z_1))$ of the equation $
(\Phi_{v_2}^{{\hat z}_2 (z_1)} \circ \Phi_{v_1}^{z_1}) (x_0^0) =
\Phi_{v_1}^{{\hat z}_1 (z_1)}(x_{-1}^0)$ w.r.t. $({\hat z}_2,{\hat
z}_1) $ and
$$ {\hat z}_1 (z_1) = z_1; \; \; \; \; \; \;  {\hat z}_2 (z_1) = {\tilde z_2} (-z_1) \n{ex8}$$
(Of course, for $z_1 < 0, $ and for the left half-plane, we could
write the equalitites which are similar to (\ref{ex7}), but we
omit
 that).

 By definition, we put
$$ \psi (z_1) = \left\{
   \begin{array}{l}
    {\hat z}_2 (z_1) \;\; \; \; \; \mbox{ if } z_1<0\\
    +\infty \; \; \; \; \; \; \; \; \; \mbox{ if } z_1 = 0\\
    {\tilde z}_2 (z_1) \;\; \; \; \; \mbox{ if } z_1>0
 \end{array} \right.  : \rr \rightarrow ]0, +\infty[ $$
 and consider the triangular system
$$\left\{
   \begin{array}{l}
    \dot z_1 = f_1 (z_1,z_2)\\
    \dot z_2 = u
 \end{array} \right.  \n{ex9}$$
 with states $(z_1,z_2) \in \rtwo$ and controls $u \in \rl,$ where
 $f_1 (z_1,z_2)$ is given by
$$ f_1 (z_1, z_2) = \left\{
   \begin{array}{l}
    {z}_2^5 \sin z_2 \;\; \; \; \; \;\; \; \; \;\; \; \; \;\mbox{ if } z_2 \le 0\\
    0 \; \; \; \; \; \; \; \; \; \; \; \; \;\; \; \; \;\mbox{ if } 0< z_2 \le \psi(z_1) \\
    (z_2{-}\psi(z_1))^5 \sin (z_2 {-} \psi(z_1)) \; \; \; \; \mbox{ if }
    z_2>\psi(z_1)
 \end{array} \right.  : \rr \rightarrow ]0, +\infty[ $$
 (if $z_1{=}0,$ then $\psi(z_1){=}+\infty,$ and $f_1(0,z_2){=}0$ for all $z_2{\ge}0,$
 by definition.)
By definition, put
$$ f(z) {:=}  \left(
   \begin{array}{c}
    f_1 (z_1,z_2)\\
    0
 \end{array} \right), \; \; \; \; \;g(z) {:=}  \left(
   \begin{array}{c}
    0\\
    1
 \end{array} \right), \; \; \; \;  \phi_N(z) {=} \left( \Phi_{v_2}^{z_2} \circ \Phi_{v_1}^{z_1}\right) (x_N^0), \; \; \; \; z{=}(z_1,z_2) {\in} \rtwo, $$
and define in $]{-} \pi, \pi[ {\times} \rr \subset \rtwo_x$ the
vector fields $A(\cdot),$ and $B(\cdot)$ by $A=(\phi_0)_{\ast} f,$
$B=(\phi_0)_{\ast} g.$ Then the system
$$ \dot x =  A(x) + B(x) u \n{ex10}$$
is well-defined in $]{-} \pi, \pi[ \times \rr, $ and satisfies
conditions (A), (B1), (B2). Finally, using the maps
$\phi_1(\cdot),$ $\phi_{-1}(\cdot),$ and the above-mentioned
symmetry (see (\ref{ex6}),(\ref{ex7}),(\ref{ex8})), we can easily
extend vector fields $A(\cdot)$ and $B(\cdot)$ onto $]{-}2 \pi,
{+} 2 \pi[ \times \rr,$ then onto $]{-}3 \pi, {+} 3\pi[ \times
\rr,$ and eventually onto $\rr\times \rr$ so that (A), (B1), (B2)
hold in $\rtwo$. Of course, system (\ref{ex10}) is not globally
feedback equivalent to a triangular system  of form (\ref{rast})
in the whole $\rtwo$ (for instance, because it is globally
equivalent to (\ref{ex9}) in $]-\pi, \pi[\times \rr$). On the
other hand, system (\ref{ex10}) satisfies conditions (A), (B1),
(B2) in the whole $\rtwo$ by the construction.

Our main result is as follows.

 {\bf Theorem 3.1} {\it Assume that vector fields
$a(\cdot),$ $b(\cdot)$ and function $\beta(\cdot,\cdot)$ are of
class $C^{n{+}1}, $  and satisfy conditions (A),(B1),(B2). Then
system (\ref{r3}) is globally controllable (in the whole $\Varm$)
in any time $[t_0,T].$}

The goal of this paper is to prove theorem 3.1.

\begin{center}
{\bf 4. The reduction of the main result to a "back-stepping"
procedure. }
\end{center}

As we see from example 3.1, system (\ref{r3}) is not globally
feedback equivalent to a system from \cite{ssp_op2} in general,
which is why the technique developed in \cite{kps2},\cite{ssp_op2}
should be at least revised essentially. However, if we want to
follow this pattern, we must first pick a point $(x^{\ast},
u^{\ast}) \in \Varm \times \rl$ around which system (\ref{r3}) is
regular, and (locally!) feedback linearizable. Using conditions
(A), and (B), we easily get the existence of a point $x^{\ast} \in
\Varm$ such that
$$ \{  b(x), [a,b] (x),..., \left({\rm ad}_a^{i-1} b \right) (x)\}
\; \; \; \; \; \mbox{ is a basis of } \Delta_{i-1} (x) $$ $$
\mbox{ for all } \; \; \; \;  x \in W(x^{\ast}), \; \; i=1,...,n
\n{app14}$$ for some neighborhood $W(x^{\ast})$ of $x^{\ast}$ in
$\Varm.$

Pick any $t_1 \in ]t_0,T[.$  In order to prove theorem 3.1, it
suffices to show that we can steer any initial state $x^0$ into
$x^{\ast}$ in time $J:=[t_0,t_1],$ and that we can steer
$x^{\ast}$ into any terminal state $x^T$ in time $I:= [t_1,T]$
w.r.t.(\ref{r3}). Next, we prove the second statement only, the
proof of the first one being similar. This statement, in turn,
follows from the following theorem 4.1.

{\bf Theorem 4.1.} {\it Let $p$ be in $\{1,...,n{-}1 \}.$ Assume
that, for every $x^T{\in}{\Varm},$ there exist  a curve $t \mapsto
z(t)$ of class $C^1(I;\Varm),$ and a map
$(t,x)\mapsto\varphi(t,x){=}(\varphi_1(t,x),...,\varphi_n(t,x)){\in}\rn$
of class $C^1$  defined in some neighborhood $E{\subset} I{\times}
\Varm$ of the curve $\{(t,z(t)){\in} I{\times}\Varm \; |\; \;
t{\in} I \}$ such that:

1) for each fixed $t{\in}I,$ the map $x \mapsto
(\varphi_1(t,x),...,\varphi_n(t,x))$ defines canonical coordinates
for system (\ref{r3}) in the corresponding neighborhood $E_t:=\{ x
\in \Varm \; |\; \; (t,x) \in E\}$ of $z(t) {\in} \Varm.$

2)   $z(\cdot),$ and $\varphi(\cdot,\cdot)$ satisfy the equalities
  $$  \frac{\partial\varphi_i (t,z(t))}{\partial x} \dot z(t) =  \frac{\partial\varphi_i (t,z(t))}{\partial x}\;  a(z(t)), \; \; \; \;\;\; \; \; \; \; \;   \; i=1,\ldots,p-1, \; \; t \in
  I;$$
  $$  \frac{d\varphi_i (t,z(t))}{dt} =0, \; \; \; \; \; \; \; \; \; \; \; \; \; \; \; \; \; \; \; \; \; \; \; \; i=1,\ldots,n,\; \; t
  \in  I; \n{r15}$$
(For $p{=}1,$ (\ref{r15}) has the form $\frac{d\varphi_i
(t,z(t))}{dt} =0, \; $ $i=1,\ldots,n,$  $t \in I,$ by definition)

3) $z(t_1)=x^{\ast},$ $z(T) {=} x^T,$ and $\frac{\partial
\varphi_p}{\partial x} (t_1, x^{\ast}) \dot z (t_1) =
\frac{\partial \varphi_p}{\partial x} (t_1, x^{\ast})
a(x^{\ast}).$

Then, for every $x^T \in \Varm,$ there exist a curve $t\mapsto
y(t)$ of class $C^1(I,\Varm),$ and a map
$(t,x)\mapsto\psi(t,x){=}(\psi_1(t,x),...,\psi_n(t,x)){\in}\rn$ of
class $C^1$ defined in some neighborhood $G{\subset} I{\times}
\Varm$ of the curve $\{(t,y(t)){\in} I{\times}\Varm \; |\; \;
t{\in} I \}$ such that

4) For each fixed $t{\in}I,$ the map
$x{\mapsto}(\psi_1(t,x),...,\psi_n(t,x))$ defines canonical
coordinates for system (\ref{r3}) in the neighborhood $G_t:=\{ x
\in \Varm \; |\; \; (t,x) \in G\}$ of $y(t) \in \Varm.$

5)  $\; y(\cdot)$ and $\psi(\cdot,\cdot)$ satisfy the equalities

 $$   \frac{\partial\psi_i (t,y(t))}{\partial x} \dot y(t) =  \frac{\partial\psi_i (t,y(t))}{\partial x}\;  a(y(t)), \;\;  \;\; \; \; \; \; \; \; \; \; \; \; \; \; i=1,\ldots,p, \; \; t \in I;
 $$
   $$ \frac{d\psi_i (t,y(t))}{dt} =0, \;  \; \; \; \; \; \; \; \; \; \; \; \; \; \; \; \; \; \; \; \; \;\; i=1,\ldots,n,\; \; t \in
   I; \n{r16}$$

6) $y(t_1)=x^{\ast},$ $y(T) {=} x^T,$ and
 $\frac{\partial\psi_{p{+}1}}{\partial x} (t_1, x^{\ast}) \dot y (t_1) =
\frac{\partial \psi_{p{+}1}}{\partial x} (t_1, x^{\ast})
a(x^{\ast})$}

Let us show that theorem 4.1 implies theorem 3.1. Indeed, for
$p=1,$ the construction of $z(\cdot),$ and $\varphi(\cdot,\cdot)$
such that conditions 1),2),3) of theorem 4.1 hold is
straightforward. Let $\zeta{=}{\tilde\varphi}(x)$
($\zeta_i{=}{\tilde\varphi_i}(x),$ $i=1,...,n$) be any canonical
coordinates for system (\ref{r3}) in a neighborhood of $x^{\ast},$
and let (\ref{rast}) be the dynamics of (\ref{r3}) in the
coordinates $(\zeta_1,...,\zeta_n).$ Given an arbitrary
$x^T{\in}\Varm,$ pick any $z(\cdot){\in} C^1(I;\Varm)$ such that
$z(t_1)=x^{\ast},$ $z(T)=x^T.$ Since the definition of the Lie
derivative does not depend on local coordinates, the condition
 $\frac{\partial {\tilde \varphi_{1}}}{\partial x} (t_1, x^{\ast}) \dot z(t_1) =
\frac{\partial {\tilde \varphi_{1}}}{\partial x} (t_1, x^{\ast})
a(x^{\ast})$ is equivalent to the equality
$\dot{\zeta_1}(t_1)=f_1({\tilde\zeta}_1^{\ast},{\tilde
\zeta}_2^{\ast})$ (where $\zeta(t)={\tilde \varphi}(z(t)),$
$\zeta_i (t) := {\tilde \varphi}_i (z(t)),$
${\tilde\zeta}_i^{\ast} = {\tilde\varphi}_i (x^{\ast})$,
$i=1,...,n$), and does not depend on the choice of canonical
coordinates around $x^\ast.$ Therefore, we can easily choose
$z(\cdot) \in C^1(I;\Varm)$ such that the condition
$\frac{\partial {\tilde \varphi_{1}}}{\partial x} (t_1, x^{\ast})
\dot z(t_1) = \frac{\partial {\tilde \varphi}_{1}}{\partial x}
(t_1, x^{\ast}) a(x^{\ast})$ holds along with the equalities
$z(t_1){=}x^{\ast},$ $z(T){=}x^T$ from the very beginning. Then
$z(\cdot)$ satisfies condition 3) with $p{=}1$ for every function
$(t,x)\mapsto\varphi(t,x) \in \rn$  defined in a neighborhood of
$\{(t,z(t))\; | \; \; t{\in} I\},$ and satisfying conditions 1),2)
of theorem 4.1. In order to construct $\varphi(\cdot,\cdot),$
consider vector fields $v_1(\cdot),...,v_n(\cdot)$ on $\Varm$ such
that $\Delta_i(x)={\rm span}
\{v_{n{-}i}(x),v_{n{-}i{+}1}(x),...,v_n(x)\}$ for all
$x{\in}\Varm,$ $i{=}0,...,n{-}1.$ (Since $\Varm$ is simply
connected, it is orientable as well as $\Delta_i(\cdot),$
$i=0,...,n{-}1,$ and such vector fields do exist). Then, for each
$\xi{\in}\Varm,$ the map $(t_1,...,t_n) \mapsto \left(
\Phi_{v_n}^{t_n} \circ \Phi_{v_{n{-}1}}^{t_{n{-}1}}
\circ...\circ\Phi_{v_1}^{t_1}\right) (\xi)$ is a diffeomorphism of
some (small) neighborhoods $\B_{\xi}(0)$ and $U(\xi)$ of
$0{\in}\rn$ and $\xi{\in}\Varm$ respectively. Let
$x\mapsto\phi(\xi,x)$ ($t_i=\phi_i(\xi,x),$ $i{=}1,...,n$) be the
inverse diffeomorphism of $U(\xi)$ onto $\B_{\xi} (0).$  For any
fixed $(t_1^0,...,t_i^0),$ the map
 $(t_{i{+}1},...,t_n) \mapsto \left( \Phi_{v_n}^{t_n} \circ
\Phi_{v_{n{-}1}}^{t_{n{-}1}} \circ ...\circ
\Phi_{v_{i{+}1}}^{t_{i{+}1}} \circ \Phi_{v_i}^{t_i^0} \circ...
 \circ \Phi_{v_1}^{t_1^0} \right) (\xi)$ defines the integral
manifold of the distribution $\Delta_{n{-}i{-}1}(\cdot)$ in
$U(\xi),$ $i{=}1,...,n{-}1.$ Therefore, for every fixed
$\xi{\in}{\Varm},$ the map $x\mapsto\phi(\xi,x)$
($t_i=\phi_i(\xi,x),$ $i{=}1,...,n$) defines canonical coordinates
$(t_1,...,t_n)$ for system (\ref{r3}) in some neighborhood
$U(\xi)$ of $\xi.$ Taking into account that $\phi(\xi,\xi)=0$ for
all $x{\in} \Varm,$ we obtain that the map $x\mapsto\varphi(t,x)$
defined by $\varphi(t,x){=}\phi(z(t),x)$ in some neighborhood of
$\{(t,z(t)) \; | \; \; t{\in}I \}$ and the curve $z(\cdot) \in
C^1(I;\Varm)$ satisfy conditions 1), 2), 3) of theorem 4.1 with
$p{=}1.$

Then, using theorem 4.1, and induction over $p=1,2,...,n{-}1,$ we
get the existence of a curve $y(\cdot) {\in} C^1(I;\Varm)$ and a
map $x{\mapsto} \psi(t,x){=}(\psi_1(t,x),...,\psi_n(t,x)) \in \rn$
of class $C^1$ in a neighborhood $G\subset I \times \Varm$ of
$\{(t,y(t)) \; | \; \; t{\in}I \}$ such that conditions 4), 5), 6)
of theorem 4.1 hold with $p{=}n{-}1,$ which implies, in particular
that $y(t_1){=}x^{\ast},$ $y(T){=}x^T,$ and
$$\frac{\partial\psi_i (t,y(t))}{\partial x} \dot y(t) =
\frac{\partial\psi_i (t,y(t))}{\partial x} a(y(t)), \;\; \;
i=1,...,n{-}1, \; \; t \in I; \n{rast2}$$ Since
$(\psi_1(t,\cdot),...,\psi_n(t,\cdot))$ are canonical coordinates
in some neighborhood of $y(t)$ for each $t{\in}I,$ we obtain from
condition (A):
$$\frac{\partial\psi_i (t,y(t))}{\partial x} b(y(t)) = 0, \;\; \;\;\; \;\; \;
\;i=1,\ldots,n{-}1, $$
$$  \frac{\partial\psi_n (t,y(t))}{\partial
x} b(y(t)) \not= 0, \;  \;\; \; \; \;  \mbox{ for all } \; \;
t{\in} I \n{rast3}$$ Therefore, (\ref{rast2}) implies that
$$ \frac{\partial\psi_i (t,y(t))}{\partial x} \dot y(t) =
\frac{\partial\psi_i (t,y(t))}{\partial x} \left(a(y(t)) +
\beta(y(t),u(t)) b(y(t)) \right), \; \; \;\; t \in I$$
 $$ \; \;
i=1,...,n{-}1, \;\; \; \; \; \; \; \; \;  \mbox{ for each } \; \;
u(\cdot) \in L_{\infty} (I; \rl)
$$
Then using a modification of the well-known Filippov lemma,
condition (A), and (\ref{rast3}) (a similar argument can be found
in \cite{pavl1}), we get the existence of $u(\cdot) {\in}
L_{\infty} (I; \rl)$ such that
$$ \frac{\partial\psi_n (t,y(t))}{\partial x} \dot y(t) =
\frac{\partial\psi_n (t,y(t))}{\partial x} (a(y(t)) {+}
\beta(y(t),u(t)) b(y(t))),\; \;  \;\;\mbox { a. e. on } \; I, $$
which yields
$$ \frac{\partial\psi (t,y(t))}{\partial x} \dot y(t) =
\frac{\partial\psi(t,y(t))}{\partial x} (a(y(t)) {+}
\beta(y(t),u(t)) b(y(t))),\; \; \;\;\mbox { a. e. on } \; I $$
Since the Jakoby matrix $\frac{\partial\psi(t,y(t))}{\partial x}$
is invertible for every $t{\in}I,$ the last equality is equivalent
to
$$\dot y (t)= a(y(t)) + \beta(y(t),u(t)) b(y(t)),\; \; \; \;\;\mbox { a. e.
on } \; \; I.$$ Therefore, $u(\cdot)$ steers $x^{\ast}$ into $x^T$
in time $I = [t_1,T]$ w.r.t (\ref{r3}). The proof of the fact that
$x^0$ can be steered into $x^{\ast}$ in time $J =[t_0,t_1]$ is
similar. Thus our main goal is to prove theorem 4.1.

\vskip10mm
\begin{center}
{\bf 5. Proof of theorem 4.1.}
\end{center}

Next we always treat $\Varm$ as $\rn,$ and always write $\rn$
instead of $\Varm$ just to simplify the notation, and to make the
arguments clearer.  For every $\varepsilon>0,$ every $\zeta{\in}
\rp,$ and every $x{\in}\rn,$ we denote by $B_{\varepsilon}
(\zeta),$ and $\Omega_{\varepsilon}(x)$ the open balls in $\rp$
and in $\rn$ respectively given by $B_{\varepsilon}
(\zeta):=\{\overline{\zeta}{\in}\rp \; |\; \;
|\overline{\zeta}{-}\zeta|<\varepsilon \};$
$\Omega_{\varepsilon}(x):=\{\overline{x}{\in}\rn \; |\; \;
|\overline{x}{-}x|<\varepsilon \},$ where $|\cdot|$ is the
standard norm generated by the standard scalar products of $\rp$
and $\rn$ respectively.

Take an arbitrary $p{\in} \{1,...,n{-}1\},$ and an arbitrary $x^T
{\in} \rn.$ Assume that a curve $z(\cdot) {\in} C^1 (I;\rn),$ and
a map $(t,x) \mapsto \varphi (t,x){=}
(\varphi_1(t,x),...,\varphi_n (t,x)) \in \rn,$ which is of class
$C^1$ in some neighborhood $E{\subset} I{\times}\rn$ of the set
$\{ (t,z(t)){\in} I{\times}\rn \; | \; \; t{\in}I \},$ satisfy
conditions 1),2),3) of theorem 4.1. Let $\zeta{=}\phi (x)$
($\zeta_i{=}\phi_i (x),$ $i{=}1,...,n$) be some fixed canonical
coordinates for system (\ref{r3}) in some small neighborhood
$U(x^{\ast}){\subset} W(x^{\ast})$ of $x^{\ast},$ where
$W(x^{\ast})$ is defined in (\ref{app14}) (for instance, we may
put $\phi(x):=\varphi (t_1,x)$ - see conditions 1), 3) of theorem
4.1), and let (\ref{rast}) be the dynamics of (\ref{r3}) in the
local coordinates $(\zeta_1,...,\zeta_n)$ (in the neighborhood
$U(x^{\ast}$). Choose $\overline{\sigma}>0$ ($\overline{\sigma}<
T{-}t_1$) such that $z(t){\in} U(x^{\ast})$ for all $t\in [t_1,
t_1{+}\overline{\sigma} ],$ and put by definition
$D:=\phi(U(x^{\ast}));$ $\zeta^{\ast} := \phi (x^{\ast});$
$\zeta_i^{\ast} := \phi_i (x^{\ast}),$ $i{=}1,...,n;$
 $\zeta^{\ast}(t) := \phi (z(t)),$
 $\zeta_i^{\ast}(t) := \phi_i(z(t)),$  $i{=}1,...,n,$ for all
$t\in [t_1, t_1{+}\overline{\sigma} ].$

Without loss of generality, we assume that $$D{=} \{
(\zeta_1,...,\zeta_n) {\in} \rn \; | \; \; |\zeta_k|<\sigma_k, \;
\; k=1,\ldots,n \}$$ with some $\sigma_k{>}0,$ $k{=}1,...,n,$ and
that every integral manifold of each $\Delta_i (\cdot)$
$(i{=}0,...,n{-}2)$ in $D=\phi(U(x^{\ast}))$ is equal to
$$\{ (\zeta_1,...,\zeta_n) {\in} \rn \; | \; \;
\zeta_k{=}\zeta_k^0, \; \;k=1,...,n{-}i{-}1; \; \;
|\zeta_k|<\sigma_k, \; \; k{=}n{-}i,...,n \}$$ with some
$\zeta_k^0,$ $k{=}1,...,n{-}i{-}1$ such that $|\zeta_k^0| <
\sigma_k.$

By the construction (see conditions 2), and 3) of theorem 4.1), we
have
$$ \dot \zeta_i^{\ast} (t)  = f_i ( \zeta_1^{\ast} (t),\ldots, \zeta_{i{+}1}^{\ast} (t)), \; \; i = 1,\ldots,p{-}1; \; \; t \in [t_1, t_1{+}\overline{\sigma} ]$$
$$\zeta_i^{\ast} (t_1) = \zeta_i^{\ast}, \; \; i=1,\ldots,n; \; \;  \;  \dot \zeta_p^{\ast} (t_1)  = f_p ( \zeta_1^{\ast} ,\ldots,  \zeta_p^{\ast}, \zeta_{p{+}1}^{\ast} ).$$
In addition,
$$ \frac{\partial f_p}{\partial \zeta_{p{+}1}} ( \zeta_1^{\ast} ,\ldots,  \zeta_p^{\ast}, \zeta_{p{+}1}^{\ast} ) \not = 0$$
- see (\ref{app14}). Therefore, there exist $\sigma \in ]0,
\overline{\sigma}[$ and $w(\cdot) \in C([t_1, t_1{+}\sigma];\rl)$
such that $w(t_1)= \zeta_{p{+}1}^{\ast},$ and
$$|w(t)| < \sigma_{p+1},\; \; \mbox{ i. e., }\; \; (\zeta_1^{\ast} (t),..., \zeta_{p}^{\ast} (t),w(t), \zeta_{p{+}2}^{\ast} (t),..., \zeta_{n}^{\ast} (t)) \in D $$
$$ \; \; \mbox{ for all }\;  t \in [t_1, t_1{+}\sigma]  \n{r17}$$
and
$$ \dot \zeta_p^{\ast} (t)  = f_p ( \zeta_1^{\ast} (t),..., \zeta_{p}^{\ast} (t), w(t)), \; \;\; \; \; \; \; \; \;  \; \; t \in [t_1, t_1{+}{\sigma} ] \n{r17_a}$$
For any ${\tilde \zeta}_{p{+}1} (\cdot){\in} C([t_1,t_1 +
\sigma];\rl)$ such that $$|{\tilde \zeta}_{p{+}1}(t)| <
\sigma_{p+1},\; \; \mbox{ i. e., }\; \; \{(\zeta_1,..., \zeta_{p},
{\tilde \zeta}_{p{+}1}(t), \zeta_{p{+}2}^{\ast}(t),...,
\zeta_{n}^{\ast} (t))\; | \; \; \zeta_i{\in}\rl,\; \; \;
i{=}1,...,p \}\cap D \not= \emptyset
$$ for all  $t \in [t_1, t_1{+}\sigma],$  we denote by
$t\mapsto\eta(t,{\tilde \zeta}_{p{+}1} (\cdot))$ the (maximal)
trajectory of the $p$ - dimensional control system
$$ \dot \zeta_i (t)  = f_i ( \zeta_1 (t),\ldots, \zeta_{i{+}1} (t)), \; \; \; \; i = 1,\ldots,p; \; \; t \in [t_1, t_1{+}{\sigma} ] \n{r18}$$
(with states $(\zeta_1,...,\zeta_p)$) with the control
$\zeta_{p{+}1} (\cdot){=}{\tilde \zeta}_{p{+}1} (\cdot), $ and
with the initial condition $\zeta_i (t_1) {=} \zeta_{i}^{\ast},$
$i{=}1,...,p.$

Since $|\frac{\partial f_i}{\partial \zeta_{i{+}1}}| {\not=} 0,$
$i{=}1,...,n$ in $D$ (see (\ref{app14})), the linearization of
(\ref{r18}) around $(\zeta_1^{\ast} (\cdot),..., \zeta_{p}^{\ast}
(\cdot), w(\cdot))$ given by
$$ \dot \chi_i (t)  = \sum\limits_{j{=}1}^{i{+}1} \frac{\partial f_i}{\partial \zeta_j} ( \zeta_1^{\ast} (t),..., \zeta_{i{+}1}^{\ast} (t))\chi_j(t), \; \; i = 1,...,p; \; \; t \in [t_1, t_1{+}{\sigma} ] \n{r19}$$
with states $(\chi_1,...,\chi_p) {\in} \rp $ and controls
$\chi_{p{+}1} {\in} \rl$ is completely controllable; therefore
there exist $p$ controls $w_i(\cdot) \in
C^1([t_1,t_1{+}\sigma];\rl),$ $i=1,...,p$ such that

$\qquad$

 ${\rm (C_1)}$ $\dot w_i (t_1) = w_i(t_1)=
w_i(t_1{+}\sigma)=\dot w_i (t_1{+}\sigma)=0,$ $i=1,...,p;$

${\rm (C_2)}$ Each control $w_i(\cdot)$ steers $0{\in}\rp$ into
$e_i = (0,...,0,1,0,...,0){\in}\rp$ (the $i$-th unit vector of the
standard basis in $\rp$) in time $I$ w.r.t. (\ref{r19}).

 $\qquad$

For every $\lambda = (\lambda_1,...,\lambda_p) {\in} \rp,$ we
define $w_{\lambda}(\cdot)$ as follows
$w_{\lambda}(t):=w(t)+\sum\limits_{j{=}1}^{p} \lambda_j w_j (t),$
$t\in [t_1, t_1{+}\sigma].$ Then, $w_{\lambda}(\cdot)
|_{\lambda=0} = w(\cdot);$
$\eta(t,w(\cdot))=(\zeta_1^{\ast}(t),...,\zeta_p^{\ast}(t)),$ for
all $t{\in}[t_1,t_1{+}\sigma],$ and, therefore, the trajectory
$t{\mapsto}\eta(t,w_{\lambda}(\cdot))$ is well-defined on
$[t_1,t_1{+} \sigma ]$ for all $\lambda$ from some small
neighborhood of $0{\in}\rp;$ $(\eta(t,w_{\lambda}
(\cdot)),w_{\lambda}(t),\zeta_{p{+}2}^{\ast}(t),\ldots,\zeta_n^{\ast}(t))
\in D$ for all $t{\in} [t_1,t_1{+}\sigma],$ and all $\lambda$ in
this neighborhood of $0{\in}\rp;$ and the map $\lambda \mapsto
F(\lambda)$ given by $F(\lambda):= \eta(t_1{+}\sigma,
w_{\lambda}(\cdot))$ is well-defined in this neighborhood of
$\lambda {=}0 {\in} \rp.$

Furthermore, from $({\rm C_2}),$ we get $\frac{\partial
F}{\partial \lambda} (0) = (0,...,0,1,0,....,0{)}^T \in \rp$ (the
$i$-th unit vector of the standard canonical basis in $\rp$);
therefore, there is $\varepsilon>0$ such that the map
$\lambda\mapsto F(\lambda)$ is a well-defined diffeomorphism of
$B_{\varepsilon} (0)$ onto some open neighborhood of
$\eta(t_1{+}\sigma, w(\cdot)) =
(\zeta_1^{\ast}(t_1{+}\sigma),\ldots,\zeta_p^{\ast}(t_1{+}\sigma)).$
Then there exists $\varepsilon_1>0$ such that
$$\overline{B_{\varepsilon_1} (\eta(t_1{+}\sigma, w(\cdot)))} \subset F(B_{\varepsilon} (0)) \; \; \; \; \; \; \; \; \;
\; \; \; \; \mbox{ and } \; \; \; \; \;
\overline{\Omega_{\varepsilon_1} (\zeta^{\ast} (t_1{+}\sigma))}
\subset D.$$ Without loss of generality, we may assume that
$\varepsilon{>}0, $ and $\varepsilon_1{>}0$ are small enough and
satisfy the condition: $$ |\alpha \zeta_{p{+}1}^{\ast} (t) +
(1-\alpha) w_{\lambda} (t) +\xi_{p{+}1}| < \sigma_{p{+}1},  \; \;
\; \; \; \;  \; \; \; \; \mbox{ whenever } \; \;
|\lambda|{<}\varepsilon, \; \; 0{\leq} \alpha{\leq}1, $$ $$
|\xi_{p{+}1}|{<}\varepsilon_1, \; \;\; \; \;
t{\in}[t_1,t_1{+}\sigma]
 \n{add__1}$$
 Fix any $\varepsilon_2>0$ such that
$$ \overline {\Omega_{\varepsilon_2} (z(t))} \subset E_t \; \; \; \mbox{ for all } t \in I, \n{r20} $$
where $E_t$ was defined in condition 1) of theorem 4.1, and
$$ \phi (\overline {\Omega_{\varepsilon_2} (z(t_1{+}\sigma))}) \subset  {\Omega_{\frac{\varepsilon_1}{2}}}
(\zeta^{\ast}(t_1{+}\sigma))  \n{r20aa} $$

{\bf Lemma 5.1.} {\it There are a curve $y(\cdot)\in
C^1([t_1{+}\sigma,T],\rn),$ and a map
$(t,x)\mapsto\psi(t,x){=}(\psi_1(t,x),...,\psi_n(t,x)){\in}\rn$ of
class $C^1$ defined in some neighborhood ${\tilde G}{\subset}
[t_1{+}\sigma,T]\times \rn$ of the curve $\{(t,x){\in}
[t_1{+}\sigma,T]{\times}\rn \; |\; \; x=y(t),\; \;  t{\in}
[t_1{+}\sigma,T] \}$ such that

1) For each fixed $t{\in}[t_1{+}\sigma,T],$ the map
$x{\mapsto}(\psi_1(t,x),...,\psi_n(t,x))$ defines canonical
coordinates for system (\ref{r3}) in the neighborhood ${\tilde
G}_t:=\{ x \in \rn \; |\; \; (t,x) \in {\tilde G} \}$ of $y(t) .$

2) $\; y(\cdot)$ and $\psi(\cdot,\cdot)$ satisfy the equalities

$$    \frac{\partial\psi_i (t,y(t))}{\partial x} \dot y(t) =  \frac{\partial\psi_i (t,y(t))}{\partial x}\;  a(y(t)), \;\;\; \; \; \;  \; i=1,...,p, \; \; t \in [t_1{+}\sigma,T];  $$ 
$$    \frac{d\psi_i (t,y(t))}{dt} =0, \; \; \; \; \; \; \; \; \; \;\; \; \; \; \; \; \; \;  i=1,...,n,\; \; t \in
[t_1{+}\sigma,T];$$

3) $y(T)=x^{T},$ and $y(t_1{+}\sigma) \in \overline
{\Omega_{\varepsilon_2} (z(t_1+\sigma))}. $  }

The proof of lemma 5.1 is given below in section 6. (Of course, it
is based on condition (B), and on definition 3.1).

Let us assume that lemma 5.1 is already proved. This allows us to
complete the proof of theorem 4.1 as follows. Put
$$\hat{\zeta}(t) := \phi (y(t)), \; \; \; \; \hat{\zeta}_i (t) := \phi_i (y(t)), \; \;
i=1,...,n, \; \; \; \; t\in [t_1{+}\sigma, t_1{+}\hat{\sigma}],
\n{r20a}$$ where $\hat{\sigma} {\in} ]\sigma, \overline{\sigma}]$
is such that $y(t) \in U (x^{\ast}),$ whenever $t {\in}
[t_1{+}\sigma, t_1{+}{\hat\sigma}].$ Then, we obtain from
(\ref{r20aa}) and from condition 3) of lemma 5.1
$$ \hat{\zeta}(t_1+\sigma) \in \Omega_{\frac{\varepsilon_1}{2}} (\zeta^{\ast}(t_1{+}\sigma)).    \n{r20b}$$
Using (\ref{add__1}), and the standard argument based on the
Gronwall-Bellmann lemma, and on the Brouwer fixed point theorem -
see \cite{lee_markus} (and \cite{kps2}, \cite{ssp_op2} for our
case) we get the existence of a control $\hat w (\cdot) $ of class
$C^1([t_1, t_1{+} \sigma]; \rl), $ ($|{\hat w} (t) |<
\sigma_{p{+}1},$ $t{\in}[t_1,t_1{+}\sigma], $ whereas the norm $
\parallel {\hat w}(\cdot) -w(\cdot)  {\parallel}_{L_1([t_1,t_1{+}\sigma];\rl)}$ should be small enough) such that

$\qquad$

$({\rm C_3})$ ${\hat w} (t_1) {=} \zeta_{p{+}1}^{\ast},$ $\frac{d
\hat w}{dt}  (t_1) {=} f_{p{+}1}
(\zeta_1^{\ast},...,\zeta_p^{\ast},\zeta_{p{+}1}^{\ast},\zeta_{p{+}2}^{\ast})$
${\hat w} (t_1{+}\sigma) {=} \hat{\zeta}_{p{+}1} (t_1{+}\sigma),$
$\frac{d{\hat w}}{dt} (t_1{+}\sigma) = \frac{d
\hat{\zeta}_{p{+}1}}{dt} (t_1{+}\sigma).$

$({\rm C_4})$ The map $\lambda \mapsto {\hat F} (\lambda) := \eta
(t_1{+}\sigma, {\hat w} (\cdot) {+} \sum\limits_{j=1}^{p}
\lambda_j w_j(\cdot))$ is well-defined for all $\lambda \in
B_{\varepsilon} (0),$ and $\overline{B_{\frac{\varepsilon_1}{2}}
(\eta(t_1{+}\sigma, w(\cdot)))}\subset {\hat F} (B_{\varepsilon}(0
)).$

$\qquad$

Then we obtain from condition $({\rm C_4})$ and from (\ref{r20b})
(see \cite{kps2},\cite{ssp_op2}) that there exists $\lambda^{\ast}
= (\lambda_1^{\ast},...,\lambda_p^{\ast})$ in $B_{\varepsilon}
(0)$ such that $({\hat \zeta}_1 (t_1{+}\sigma),...,{\hat \zeta}_p
(t_1{+}\sigma)) {=} \eta (t_1{+}\sigma, {\hat w}_{\lambda^{\ast}}
(\cdot)),$ where ${\hat w}_{\lambda^{\ast}} (\cdot)$ is given by
${\hat w}_{\lambda^{\ast}} (t) {:=} {\hat w} (t){+}
\sum\limits_{j=1}^{p} \lambda_j^{\ast} w_j(t),$ $t\in [t_1,
t_1{+}\sigma].$ Let us define $y(t),$ and $\overline{\zeta} (t) :=
\phi (y(t))$ on $[t_1, t_1{+}\sigma]$ as follows: by definition,
put
$$({\overline \zeta}_1 (t),...,{\overline \zeta}_p (t)) :=
\eta (t,{\hat w}_{\lambda^{\ast}} (\cdot)),\; \;  \; \; \; \; \;
\; {\overline \zeta}_{p{+}1} (t) := {\hat w}_{\lambda^{\ast}} (t),
\; \; \; \; \; \; t \in [t_1, t_1{+}\sigma]; \n{r21}
$$
in addition, let ${\overline \zeta}_{p{+}2} (\cdot),$
...,${\overline \zeta}_n (\cdot)$ be any functions of class
$C^1([t_1,t_1+\sigma];\rl)$ such that
$$ {\overline \zeta}_i (t_1{+}\sigma) = {\hat \zeta}_i (t_1{+}\sigma), \; \; \; \;
\frac{d {\overline \zeta}_i}{dt} (t_1{+}\sigma) = \frac{d {\hat
\zeta}_i}{dt} (t_1{+}\sigma), \; \; \; \; i=p{+}2,...,n, \n{r22}$$
$$  {\overline \zeta}_i (t_1) = \zeta_i^{\ast}, \; \; \; \; i=p{+}2,...,n, \n{r22a}$$
and such that
$$|{\overline \zeta}_i (t)|<\sigma_i, \; \; \; i{=} p{+}2, \ldots, n,\mbox { i. e., } \; \;
 ({\overline \zeta}_1 (t),\ldots,{\overline \zeta}_n (t)) \in D \;\; \;\mbox{ for all } t \in [t_1,t_1{+}\sigma], \n{r23}$$
 and then put:
 $$ \overline{\zeta} (t) : = (\overline{\zeta}_1 (t),...,\overline{\zeta}_n (t)),
 \;\;
\; \; \;  y(t):={\phi}^{-1} (\overline{\zeta} (t)), \; \; \; \;\;
\; \mbox{ for all } \;  t \in [t_1, t_1{+}\sigma[ \n{r24}
$$
Thus, we have constructed $y(\cdot)$ of class $C^1(I;\rl)$ such
that conditions 4), 5), 6) of theorem 4.1 hold with every
canonical coordinate functions $\psi_j (t,y),$ $j=1,...,n.$
Indeed, the inclusion $y(\cdot) \in C^1 (I;\rn) $ follows from
(\ref{r20a}), from (\ref{r24}), and from (\ref{r22}),(\ref{r21})
(in addition, we take into account that
$${\overline\zeta}_{p{+}1}(t_1{+}\sigma) = {\hat \zeta}_{p{+}1} (t_1{+}\sigma), \; \; \; \; \;\; \;
\frac{d{\overline\zeta}_{p{+}1}}{dt}(t_1{+}\sigma) = \frac{d{\hat
\zeta}_{p{+}1}}{dt} (t_1{+}\sigma)  $$ by $({\rm C_1}),$ $({\rm
C_3}),$ and by (\ref{r21}), and that
$$ \frac{d{\hat \zeta}_i}{dt} (t_1+\sigma)= f_i ({\hat
\zeta}_1 (t_1+\sigma),...,{\hat \zeta}_{i{+}1} (t_1+\sigma)) = f_i
({\overline \zeta}_1 (t_1+\sigma),...,{\overline \zeta}_{i{+}1}
(t_1+\sigma))= $$ $$ = \frac{d{\overline \zeta}_i}{dt}
(t_1+\sigma), \; \; \;\; \; \; \; \; \; \; \; \; \; \; \; \; \;
i=1,\ldots,p
$$
 by the construction). Conditions 4), and 5) of
theorem 4.1 follow from conditions 1), and 2) of lemma 5.1
respectively  (by the construction, (\ref{r16}) is true for all
$t\in [t_1, t_1{+}\sigma],$ and we can easily construct the
appropriate $\psi(t,x)$ for all $t \in I$ following the argument
from section 4). The equality $y(T) = x^T $ follows from condition
3) of lemma 5.1; the equality $y(t_1)=x^{\ast}$ follows from the
definition of $\overline{\zeta}(\cdot):$ indeed, by the definition
of $\eta(t,v(\cdot)),$ we have
$$(\overline{\zeta}_1 (t_1),...,\overline{\zeta}_p (t_1))=\eta(t_1,{\hat w}_{\lambda^{\ast}} (\cdot))=(\zeta_1^{\ast},...,\zeta_p^{\ast}); $$
conditions $({\rm C_3}),$ and $({\rm C_1})$ yield
$\overline{\zeta}_{p{+}1} (t_1) ={\hat w}_{\lambda^{\ast}} (t_1) =
\zeta_{p{+}1}^{\ast};$ taking into account (\ref{r22a}), we obtain
$\overline{\zeta}(t_1)=(\zeta_1^{\ast},...,\zeta_n^{\ast})
=\zeta^{\ast}, $ which implies that
$y(t_1)=\phi^{-1}(\overline{\zeta} (t_1)) = {\phi}^{-1}
(\zeta^{\ast}) = x^{\ast}.$

Finally, we obtain from ${\rm (C_3)},$ and from ${\rm (C_1)}$
$$\frac{d {\hat w}_{\lambda^\ast}}{dt}  (t_1) = f_{p+1}
(\zeta_1^{\ast},...,\zeta_p^{\ast},\zeta_{p{+}1}^{\ast},\zeta_{p{+}2}^{\ast}),$$
which yields: $\frac{\partial \psi_{p{+}1}}{\partial x} (t_1,
x^{\ast}) \dot y (t_1) = \frac{\partial \psi_{p{+}1}}{\partial x}
(t_1, x^{\ast}) a(x^{\ast}) $ (because the definition of the Lie
derivative does not depend on coordinates). Therefore, $y(\cdot),$
and $\psi(\cdot,\cdot)$ satisfy condition 6) of theorem 4.1 as
well. The proof of theorem 4.1 is complete.

\vskip10mm
\begin{center}
{\bf 6. Proof of lemma 5.1.}
\end{center}

Consider the following control system of ordinary differential
equations
$$ \left\{
\begin{array}{l}
    \frac{\partial\varphi_j }{\partial x}(t,x(t)) \dot x(t) =  \frac{\partial\varphi_j }{\partial x} (t,x(t))a(x(t)), \; \; \; \; \; \;\;\;\;\; j{=}1,...,p{-}1, \; \; t \in I, \; \; (t,x) \in E;  \\
    \frac{\partial \varphi_p}{\partial t} (t,x(t)) + \frac{\partial\varphi_p }{\partial x}(t,x(t)) \dot x(t) = v(t); \; \;
    \;\; \; \; \; \; \; \; \; \; \; \;\; \; \; \; \; \; \; \; \; \;\; \; \; \; \; \;\; \; \;
    t{\in}I, \; \; (t,x) \in E\\
    \frac{\partial \varphi_j}{\partial t} (t,x(t)) + \frac{\partial\varphi_j }{\partial x}(t,x(t)) \dot x(t) = 0, \; \; \; \; \; \;\; \; \; \; \; \;\; \; \;\; \;\; \;  j{=}p{+}1,...,n,\; \; t \in I, \; \; (t,x) \in E;
 \end{array}\right.  \n{rr21} $$
(If $p=1,$ then the first row is empty by definition) with states
$x{\in}\rn,$  $(t,x) {\in} E,$ and controls $v{\in}\rl.$ Since the
Jakoby matrix $\frac{\partial \varphi}{\partial x} (t,x)$ is
invertible for all $(t,x){\in} E $ (see condition 1) of theorem
4.1), we can rewrite (\ref{rr21}) in its standard form
$$\dot x (t) = \digamma (t, x(t),v(t)), \; \; \; \; \; \; \; \; \; \; \; t\in
I, \; \; (t,x) \in E, \n{rr30}$$ where $\digamma
(\cdot,\cdot,\cdot)$ is given by

$$ \digamma (t,x,v) = [\frac{\partial \varphi (t,x)}{\partial x}]^{-1}
\left( \begin{array}{c}
    \frac{\partial \varphi_1 }{\partial x}(t,x) a(x)  \\
    \ldots     \\
    \frac{\partial\varphi_{p{-}1} }{\partial x}(t,x) a(x)\\
    v - \frac{\partial \varphi_p}{\partial t} (t,x) \\
    - \frac{\partial \varphi_{p{+}1}}{\partial t} (t,x)\\
    \ldots \\
     - \frac{\partial \varphi_n}{\partial t} (t,x)
   \end{array} \right)  \n{rr31} $$

   Given $(\tau,\tilde{x}) \in E$ and $ v(\cdot) {\in} L_{\infty}
   (I;\rl),$ denote by $t{\mapsto}x(t,\tau,\tilde{x},v(\cdot))$
   the maximal trajectory of system (\ref{rr21}) with the control
   $v(\cdot)$ and with the initial condition $x(\tau,\tau,\tilde{x},v(\cdot))=\tilde{x}$
   (of course, $(t,x(t,\tau,\tilde{x},v(\cdot))) \in E$ for all
   admissible $t$).
In addition, if $v = v(t,x)$ is a feedback control, which can be
time-varying, and even discontinuous, defined in some open subset
$\tilde{E} \subset E,$ in general, and if $(\tau,{\tilde{x}}) \in
\tilde{E},$ then we denote by $t\mapsto x(t,\tau,\tilde{x},
v(\cdot,\cdot))$ the (maximal) trajectory of (\ref{rr21}) such
that $x(\tau,\tau,\tilde{x}, v(\cdot,\cdot)) = \tilde{x}$ as well
(if this trajectory is well-defined).  From conditions 2) and 3)
of theorem 4.1, it follows that
   $$z(t) = x(t,T, x^T, v_0(\cdot)) \; \; \; \; \mbox { for all }  \; \; t\in I, \n{rr22}$$
  with
   $$ v_0(t):=0 , \; \; \; \; \; \; \; \; t\in I, \n{rr23}$$
   and then, using the Gronwall-Bellmann lemma, we get the
   existence of $\delta>0$ such that the trajectory $t\mapsto x(t,T,x^T,v(\cdot))$
   of system (\ref{rr30}) is well-defined for all $t{\in}I,$ and
   $$ (t,x(t,T,x^T,v(\cdot))) \in E \; \; \; \;  \mbox{ and } \; \; \; \; |x(t,T,x^T, v(\cdot))-z(t)| < \frac{\varepsilon_2}{2}  $$
   $$   \mbox{ for all } t{\in}I, \; \; \; \; \; \; \; \; \; \; \;  \mbox{ whenever } \;  \parallel v(\cdot) - v_0 (\cdot) {\parallel}_{L_{\infty}(I;\rl) } < \delta  \n{rr24}$$
for every $v(\cdot) \in L_{\infty} (I; \rl)$ where
$\varepsilon_2>0$ was defined in (\ref{r20aa}).

By definition, put:
 $$\T := \{ (t,x) \in I\times\rn \; | \; \; \; t {\in}I; \; \; |x{-}z(t)|\le
 \varepsilon_2
 \}. \n{L1add}$$

{\bf Lemma 6.1.} {\it There exist a curve $x(\cdot) {\in} C (I;
\rn),$ a finite sequence of open sets $\{\T_i{\}}_{i=1}^{N}$ of
the form $\T_i
=]\tau_i{-}\alpha_i,\tau_i{+}\alpha_i[\times\Omega_{\beta_i}(x_i)$
with some $(\tau_i,x_i) \in \T,$ $i=1,...,N,$ a finite sequence of
numbers $\tau_i^{\ast}$ in $[t_1,T],$ $i=1,\ldots,N{+}1,$ $N$
finite sequences of vector fields $\nu_i {=} \{ \nu_i^k
(\cdot){\}}_{k=1}^{k_i}$ each of which belongs to $\Delta_{n-p-1}
(\cdot)$ (i.e. $\nu_i^k(\cdot) {\in} \Delta_{n{-}p{-}1}(\cdot),$
$k{=}1,...,k_i,$ $i{=}1,...,N$), and $N$ finite sequences of
numbers $\mu_i {=} \{ \mu_i^k {\}}_{k=1}^{k_i},$  $\mu_i^k {\ge} 0
,$ $k{=}1,...,k_i,$ $i{=}1,...,N,$
 such that

1) $\overline{\Omega_{\beta_i} (x_i)} \subset \D_{{\nu}_i^1},\; $
$k=1,...,k_i;\;\; \; \;  $ $ \bigcup\limits_{i=1}^{N} \T_i \subset
E,\; $ $i{=}1,...,N$

2) $ \; \tau_1^{\ast}{=}T{>} {\tau}_2^{\ast}
{>}...{>}{\tau}_N^{\ast}{>}{\tau}_{N{+}1}^{\ast}{=}t_1; \; \; $
$x(T){=}x^T; \; $ and $(t,x(t))\in \T_i,$ whenever
$t{\in}[\tau_{i{+}1}^{\ast},\tau_{i}^{\ast}],$ $i{=}1,...,N$

3) For every $i{=}1,...,N,$ and every $t\in ]\tau_{i{+}1}^{\ast},
\tau_{i}^{\ast}[, \;\;  $ $\dot x (t)$ is well-defined, and

$$ \frac{\partial \varphi_j}{\partial x} (t,x(t)) \dot x (t) = \frac{\partial \varphi_j}{\partial x} (t,x(t)) a(x(t)), \; \; \; \; \; \; j{=}1,...,p{-}1, \; \; \; t\in ]\tau_{i{+}1}^{\ast}, \tau_i^{\ast}[, \; \; \; i{=}1,...,N, \n{L1}$$
$$ \frac{\partial \varphi_p}{\partial t} (t,x(t)) +  \frac{\partial \varphi_p}{\partial x} (t,x(t)) \dot x (t) = v_i(t,x(t)), \; \; \; \; \; \;\; \; \; \; \; \; \; \;\; \; \; \;  \; \; \; t\in ]\tau_{i{+}1}^{\ast}, \tau_i^{\ast}[, \; \; \; i{=}1,...,N, \n{L2}$$
$$ \frac{\partial \varphi_j}{\partial t} (t,x(t)) +  \frac{\partial \varphi_j}{\partial x} (t,x(t)) \dot x (t)= 0, \; \; \; \; \; \;\; \; \; \; \; \; j{=}p{+}1,...,n, \; \; \; t\in ]\tau_{i{+}1}^{\ast},\tau_i^{\ast}[, \; \; \; i=1,...,N, \n{L1_a_a}$$

where $v_i(t,x),$ $i=1,...,N$ are given by

$$ v_i(t,x) = \frac{\partial \varphi_p}{\partial t} (t,x ) + \frac{\partial \varphi_p}{\partial x} (t,x ) \left(\Phi_{\nu_i}^{\mu_i}\right)_{\ast} a(\Phi_{\nu_i}^{-\mu_i}(x))
\; \; \; \; \; \; \; \mbox{ for all }\; \;  (t,x) {\in} \T_i,  \;
\; i{=}1,...,N \n{L3}
$$
and satisfy the conditions
$$|v_i (t,x)| < \delta \; \; \; \;  \; \; \; \; \; \; \; \;\; \; \; \;  \mbox{ for each }\; \; \;  (t,x) \in \T_i, \; \; \; \; \; i{=}1,...,N \n{L4}$$
with $\delta$ defined in (\ref{rr24}).
 }

(According to the definition of the diffeomorphisms
$\Phi_{\nu_i}^{\mu_i} (\cdot) $ and $\Phi_{\nu_i}^{-\mu_i}
(\cdot)$ for finite sequences $\nu_i {=} \{\nu_i^k (\cdot)
{\}}_{k=1}^{k_i},$ $\mu_i {=} \{\mu_i^k {\}}_{k=1}^{k_i}$ which
was given in section 2, we have: $\Phi_{\nu_i}^{-\mu_i}  {=}
\Phi_{\nu_i^{k_i}}^{-\mu_i^{k_i}}{\circ} ...{\circ}
\Phi_{\nu_i^1}^{-\mu_i^1} ;\; $   $\Phi_{\nu_i}^{\mu_i} {=}
\Phi_{\nu_i^{1}}^{\mu_i^{1}} {\circ} ... {\circ}
\Phi_{\nu_i^{k_i}}^{\mu_i^{k_i}}$)

First we assume that lemma 6.1 is already proved, and prove lemma
5.1. The proof of lemma 6.1, in turn, is based on condition (B)
for system (\ref{r3}) and is given in Appendix (we point out that
it is a modification of the proofs of lemmas 3.4, and 3.1.1  from
\cite{ssp_op2}; in particular, condition (B) allows us to find
$v_i(t,x)$ given by (\ref{L3}), and satisfying (\ref{L4})). Let
$x(\cdot)$ be a curve from lemma 6.1. To make the proof of lemma
5.1 clearer, we assume that $k_i=1,$ $i=1,...,N,$ i.e. each
sequence $\mu_i = \{ \mu_i^k {\}}_{k=1}^{k_i}$ and $\nu_i =\{
\nu_i^k (\cdot) {\}}_{k=1}^{k_i}$ consists of one element only,
and then, to simplify the notation, we put $\mu_i:=\mu_i^1,$
$\nu_i(\cdot):=\nu_i^1(\cdot)$ (however, we will explain how we
can adjust our construction in the general case)

Put $v(t):= v_i (t,x(t)),$ $t {\in} ]\tau_{i{+}1}^{\ast},
\tau_{i}^{\ast}],$ $i=1,...,N.$ Then, $v(\cdot)$ is a piecewise
continuous open-loop control. Combining (\ref{L4}), and
(\ref{rr24}), and taking into account (\ref{L1})-(\ref{L1_a_a}),
and condition 2) of lemma 6.1, we get
$$ |x(t)-z(t)|< \frac{\varepsilon_2}{2} \; \; \; \;  \mbox{ for all } t \in I \n{L4a}  $$

Let $N_0$ in $\{1,...,N \}$ be such that $\tau_{N_0{+}1}^{\ast}
{\leq} t_1{+}\sigma {<} \tau_{N_0}^{\ast}.$ Without loss of
generality, we assume that $\tau_{N_0{+}1}^{\ast}{=}t_1{+}\sigma;$
otherwise, with slight abuse of notation, we put by definition
$\tau_{N_0{+}1}^{\ast}{:=}t_1{+}\sigma,$ whereas
$\tau_{i}^{\ast},$ $i{=}1,...,N_0$ are the same (the terminal
point of the curve $y(\cdot)$ mentioned in lemma 5.1 is
$t_1{+}\sigma,$ and, therefore, we should deal with
$[t_1{+}\sigma,T]$ instead of $[t_1,T]$ in this section).

Take any $\varepsilon_3>0$ such that $\varepsilon_3<
\frac{\varepsilon_2}{2}.$  For any $\tau{\le}t$ in
$[t_1{+}\sigma,T],$ we put by definition:
$$ \Gamma_{t}^{\tau} = \{ (s,z) \in [\tau,t] \times \rn\; |\; \; s{\in}[\tau,t]; \; \; |z-x(t)|<\varepsilon_3  \}; \; \; \; \; \; \; \; \; \Gamma:=\Gamma_{t_1{+}\sigma}^T  \n{L5}$$
By the construction, $(\tau_{i{+}1}^{\ast},
x(\tau_{i{+}1}^{\ast})) \in \T_i \cap \T_{i+1},\; $
$i{=}1,...,N_0{-}1;$ therefore $\varepsilon_3$ in
$]0,\frac{\varepsilon_2}{2}[$ can be chosen such that $\Gamma$ in
(\ref{L5}) satisfies the conditions
$${\overline \Gamma}_{\tau_{1}^{\ast}}^{\tau_{1}^{\ast}} \subset \T_1; \; \; \; \; {\overline \Gamma}_{\tau_{i+1}^{\ast}}^{\tau_{i+1}^{\ast}} \subset \T_i{\cap}\T_{i{+}1}, \; \; i=1,...,N_0-1; \; \; \; \;{\overline \Gamma}_{\tau_{N_0+1}^{\ast}}^{\tau_{N_0+1}^{\ast}} \subset \T_{N_0}; \; \; \; \;
\n{L6}$$
$$ \Gamma \subset \bigcup\limits_{i=1}^{N_0} \T_i. \n{L7}$$

Put $\sigma^0:= \min \{
\frac{\tau_i^{\ast}-{\tau_{i{+}1}^{\ast}}}{2}, \; \;1 {\le}  i
{\le} N_0 \}.$ Then, there exists $M{>}0$ such that, for every
sequence $\varkappa {=} \{\sigma_i {\}}_{i{=}1}^{N_0{+}1}$
satisfying the conditions $0{<} \sigma_i {<}\sigma^0, $
$i{=}1,...,N_0{+}1,$ there are a smooth time-varying vector field
$\nu_{\varkappa} (t,\cdot),\; $ $t{\in}[t_1{+}\sigma,T],$ and a
smooth function $\mu_{\varkappa} (\cdot) \in C^{\infty}
([t_1{+}\sigma;T];\;  [0,+\infty[)$ such that
$$ \mu_{\varkappa} (T) {=} \mu_{\varkappa}(t_1{+}\sigma) {=} 0; \; \; \;  \; \;\; \mbox{ and } \;  \; \;
\mu_{\varkappa}(t)=\mu_i,$$  $$  \mbox{ whenever } \; \;  \; \; t
\in [\tau_{i{+}1}^{\ast}{+}\sigma_{i{+}1}, \;
\tau_{i}^{\ast}{-}\sigma_{i}], \; \; i=1,\ldots, N_0;  \n{L8}
$$
$$ \nu_{\varkappa} (t, \cdot) \in \Delta_{n{-}p{-}1} (\cdot),\; \;
\; \mbox{ whenever } t \in [t_1{+}\sigma,T]; \n{L9}$$
$$\nu_{\varkappa} (t, x) = \nu_i(x), \;\; \; \; \; \; \mbox{ whenever } \; \; x{\in}\D_{\nu_i},\; \;
 t \in [\tau_{i{+}1}^{\ast}{+}\sigma_{i{+}1}, \;
\tau_{i}^{\ast}{-}\sigma_{i}], \; \; i=1,\ldots,N_0,  \n{L10}$$
and such that the feedback control $v_{\varkappa}(t,x)$ given by
$$ v_{\varkappa}(t,x) = \frac{\partial \varphi_p}{\partial t} (t,x ) + \frac{\partial \varphi_p}{\partial x} (t,x ) \left(\Phi_{\nu_{\varkappa}(t) }^{\mu_{\varkappa}(t)}\right)_{\ast} a(\Phi_{\nu_{\varkappa}(t)}^{-\mu_{\varkappa}(t)}(x))
\n{L11}$$ satisfies the condition
$$ \max \{ |\digamma (t,x,v_{\varkappa}(t,x))|\;  \; |\; \; t{\in} [t_1{+}\sigma;T]; \; \; x{\in}\overline{\Omega_{\varepsilon_3} (x(t))}  \} \leq M, \n{L12}$$
where $M $ is given by
$$ M:= \max\{ |\digamma (t,x,v)|\; \; |  \; \; \; v =  \frac{\partial \varphi_p}{\partial t} (t,x ) + $$
$$ +  \frac{\partial \varphi_p}{\partial x} (t,x )
\left(\Phi_{\nu_{i}}^{\mu}\right)_{\ast}
a(\Phi_{\nu_{i}}^{-\mu}(x)), \; \; \; \; (t,x) \in
\overline{\T_i}, \; \; 0 {\le} \mu {\le} \mu_i; \; \;
i{=}1,...,N_0 \}. \n{L12a}
$$
For instance, given $\varkappa {=} \{\sigma_i
{\}}_{i{=}1}^{N_0{+}1}$ with small enough
$\sigma_i{\in}]0,\sigma^0[,$ take any functions $\lambda_i(\cdot)
\in C^{\infty} ([t_1{+}\sigma,T];\rr),$ $i=1,...,N_0,$ and
$\lambda(\cdot) \in C^{\infty} ([t_1{+}\sigma,T];\rr)$ such that
$$ \sum\limits_{i=1}^{N_0} \lambda_i (t) =1, \; \; \; \;  \mbox{ and } \; \; \; \;  \lambda_i(t)\ge 0, \; \; \;i=1,...,N_0, \; \; \; \;  \; \; \; \mbox{ for all }\;  t\in [t_1{+}\sigma,T]; $$
$$\lambda_i(t) = 1, \; \; \; \; \; \; \; \mbox{ whenever }\;  t\in  [\tau_{i{+}1}^{\ast} {+} \frac{\sigma_{i{+}1}}{2}, \;  \tau_{i}^{\ast}
{-} \frac{\sigma_{i}}{2}], \; \; \; i=2,\ldots, N_0-1;
 \; \; \;  \lambda_1 (t) = 1, $$ $$  \mbox{ whenever } \; t\in  [\tau_{2}^{\ast} {+} \frac{\sigma_{2}}{2}, \;
\tau_{1}^{\ast}];   \; \; \; \; \; \; \; \lambda_{N_0} (t) = 1 \;
\; \; \; \; \;  \mbox{ whenever } \;  t\in [\tau_{N_0{+}1}^{\ast},
\; \tau_{N_0}^{\ast} {-} \frac{\sigma_{N_0}}{2}];
 $$
$$\lambda_i(t) +  \lambda_{i+1}(t) = 1, \; \; \; \; \; \; \; \mbox{ whenever } t\in  [\tau_{i{+}1}^{\ast} {-} \frac{\sigma_{i{+}1}}{2}, \; \tau_{i{+}1}^{\ast}
{+} \frac{\sigma_{i{+}1}}{2}], \; \; \; i=1,\ldots,N_0{-}1;
 $$
$$\lambda(t) = 1 \; \; \; \; \; \; \; \; \; \mbox{ whenever } \; t\in  \bigcup\limits_{i{=}1}^{N_0}[\tau_{i{+}1}^{\ast} {+} \sigma_{i{+}1}, \; \tau_{i}^{\ast}
{-} {\sigma_{i}}] ;$$
$$ 0\le \lambda(t) \le 1 \; \;\; \; \; \; \; \; \; \; \; \; \;  \; \mbox{ whenever } \;\;   t \in [t_1{+}\sigma,T]; $$
$$ \lambda(t) = 0 \; \; \; \; \; \; \; \; \; \mbox{ whenever } \;  \; t \in
\bigcup\limits_{i=1}^{N_0{-}1}[\tau_{i{+}1}^{\ast} {-}
\frac{\sigma_{i{+}1}}{2}, \; \tau_{i{+}1}^{\ast} {+}
\frac{\sigma_{i{+}1}}{2}];$$
$$\lambda(T)=\lambda(t_1{+}\sigma) = 0. $$
Then $\mu_{\varkappa}(\cdot)$ and $\nu_{\varkappa} (\cdot,\cdot)$
given by $\mu_{\varkappa}(t) {=} \sum\limits_{i=1}^{N_0}\mu_i
\lambda_i (t) \lambda(t),$ and $\nu_{\varkappa}(t,x) {=}
\sum\limits_{i=1}^{N_0} \lambda_i (t) \nu_i(x),$ $t {\in}
[t_1{+}\sigma,T],\; $ $x\in\D_{\nu_i}$ satisfy
(\ref{L8})-(\ref{L12a}).

Let us remark that our vector field $\nu_{\varkappa} (t)$ is
actually time-varying only around the moments of switching
$\tau_i^{\ast},$ which allows us to construct the smooth feedback
control (\ref{L11}). If each $\nu_i$ were a sequence of vector
fields $\{\nu_i^k (\cdot){\}}_{k=1}^{k_i},$ we would have to take
into account each switching from $\nu_i^k (\cdot)$ to
$\nu_i^{k{+}1} (\cdot).$ The above-mentioned convex combinations
would become more complicated, but the idea would be the same.

To simplify the notation, put
$\nu_{\varkappa}(t):=\nu_{\varkappa}(t,\cdot).$

{\bf Lemma 6.2} {\it There exists $\varkappa{=}\{\sigma_i
{\}}_{i=1}^{N_0{+}1}$ with small enough $\sigma_i,$
$0{<}{\sigma_i}{<}\sigma^0,$   $i{=}1,...,N_0{+}1$ such that the
corresponding trajectory $t\mapsto x_{\varkappa}(t):=x(t,T,x^T,
v_{\varkappa}(\cdot,\cdot))$ of system (\ref{rr21}) with the
(smooth) feedback control $v_{\varkappa}(\cdot,\cdot)$ given by
(\ref{L11}) and with the "initial" condition $x_{\varkappa}(T)=
x^T$ is well-defined for all $t \in [t_1{+}\sigma, T],$ and
satisfies the condition
$$|x_{\varkappa} (t) - x(t) | < \varepsilon_3 \; \; \; \mbox{ for all } \; t \in [{t_1}{+}\sigma,T] \n{L12_b}$$
(which implies $(t,x_{\varkappa}(t))) \in \Gamma,$ whenever $t\in
[t_1{+}\sigma, T]$). }

To prove lemma 6.2, we just note that
$$x_{\varkappa} (t) {=} x(t,T,x^T,v_{\varkappa}(\cdot,\cdot)) \; \; \; \;\mbox{ and }  \; \; \; x(t){=}x(t,T,x^T,v(\cdot,\cdot)), \; \;
\; \; \mbox{ for all } \;  \; t {\in} [t_1{+}\sigma,T] $$ and
$v_{\varkappa}(t,x)=v(t,x),$ whenever $t {\in}
\bigcup\limits_{i=1}^{N_0} [\tau_{i{+}1}^{\ast}{+} \sigma_{i{+}1},
\; \tau_{i}^{\ast} {-} \sigma_{i}],\; $ $x \in
\Omega_{\varepsilon_3} (x(t)),$ where
$$ v (t,x) := v_i (t,x), \; \; \; \; \; \; \; \; \mbox{ whenever } \; \tau_{i{+}1}^{\ast}{<}t {\le} \tau_i^{\ast}, \; \; (t,x){\in}\T_i, \; \;  i=1,\ldots,N_0, \n{vstavka}$$
with $v_i(t,x)$ defined in (\ref{L3}). In addition,
$$ \max \{ |\digamma (t,x,v_{\varkappa}(t,x))| \; \;| \;  \; \;  t{\in}[t_1{+}\sigma,T], \; \; x\in\overline{\Omega_{\varepsilon_3} (x(t))}   \} \leq M; $$
$$ \max \{ |\digamma (t,x,v(t,x))| \; \;| \;  \; \;  t{\in}[t_1{+}\sigma,T], \; \; x\in\overline{\Omega_{\varepsilon_3} (x(t))}   \} \leq M; $$
Therefore, if $\sigma_i>0$ are small enough, then
$t{\mapsto}x_{\varkappa}(t)$ is well-defined on
$[t_1{+}\sigma,T],$ and $\parallel x_{\varkappa} (\cdot) -x(\cdot)
{\parallel}_{C(I;\rn)}$ is small enough, which can be proved by
the standard argument based on the Gronwall-Bellmann lemma. The
proof of lemma 6.2 is complete.

Finally, we put:  $\mu(t){:=}\mu_{\varkappa} (t),$
$\nu(t){:=}\nu_{\varkappa} (t),$ $t{\in}[t_1{+}\sigma,T],$ with
$\mu_{\varkappa} (t),$  $\nu_{\varkappa} (t)$ from lemma 6.2, and

$$ y(t): = \Phi_{\nu(t)}^{-\mu(t)} (x_{\varkappa}(t)), \; \; \; \; \;\; \; \; \; \; \;  \; \mbox { whenever } \; \; t \in [t_1{+}\sigma,T] \n{L13} $$
$$ \tilde{\psi}_j (t,y):= \varphi_j (t, \Phi_{\nu(t)}^{\mu(t)} (y)), \; \;  \; \; \; \; \; \; j=1,\ldots,p,\; \; \; \;\; \; \; \mbox{ whenever } \;
\; y \in \Phi_{\nu(t)}^{-\mu(t)} (
\Omega_{\varepsilon_3}(x_{\varkappa}(t))) \n{L14}$$ Let us show
that $y(\cdot)$ defined by (\ref{L13}) satisfies conditions
1),2),3) of lemma 5.1. Indeed, taking into account that
$x_{\varkappa}(t) = \Phi_{\nu(t)}^{\mu(t)} (y(t)),$ we obtain from
(\ref{L13}), and from (\ref{rr21}):
$$\frac{\partial \varphi_j (t,\Phi_{\nu(t)}^{\mu(t)} (y(t)))}{\partial x}
\left(  \left(\frac{\partial\Phi_{\nu(t)}^{\mu(t)}}{\partial
t}\right) (y(t)) + \left( \Phi_{\nu(t)}^{\mu(t)}\right)_{\ast}
\dot y (t) \right) =
$$
$$ = \frac{\partial \varphi_j (t,\Phi_{\nu(t)}^{\mu(t)} (y(t))) }{\partial x} a(\Phi_{\nu(t)}^{\mu(t)} (y(t))),
\; \; \; \; j{=}1,...,p{-}1, \; \; t{\in} [t_1{+}\sigma,T] \n{L15}
$$
By the construction, $\nu(t) \in \Delta_{n{-}p{-}1},$ for all $t;$
therefore
$$\frac{\partial \varphi_j (t,\Phi_{\nu(t)}^{\mu(t)} (y))}{\partial x}  \left(\frac{\partial \Phi_{\nu(t)}^{\mu(t)}(y)}{\partial
t}\right) = 0 \; \; \; \; \; \;  j{=}1,...,p \n{L16}$$ for every
admissible $y$ and $t.$ In addition, from (A),(B) it follows that
$$a \left(\Phi_{\nu(t)}^{\mu(t)} (y(t))\right) - \left(\Phi_{\nu(t)}^{\mu(t)} \right)_{\ast} a(y(t)) \in \Delta_{n-p} (\Phi_{\nu(t)}^{\mu(t)} (y(t)))$$
Therefore, we get from (\ref{L15}), (\ref{L16})
$$ \frac{\partial \varphi_j (t,\Phi_{\nu(t)}^{\mu(t)} (y(t)))}{\partial x} \left( \Phi_{\nu(t)}^{\mu(t)} \right)_{\ast} \dot y (t) =
\frac{\partial \varphi_j (t,\Phi_{\nu(t)}^{\mu(t)}
(y(t)))}{\partial x}  \left( \Phi_{\nu(t)}^{\mu(t)} \right)_{\ast}
a(y(t)), \; \; \; \; j{=}1,...,p{-}1, $$
$$ t\in[t_1{+}\sigma,T].$$
Combining this with (\ref{L16}), (\ref{rr21}), (\ref{L11}),
(\ref{L13}) we obtain that the last equality holds for $j=p$ as
well. On the other hand, by the definition of $\tilde{\psi_j}
(t,y),$ the obtained equalities
$$ \frac{\partial \varphi_j (t,\Phi_{\nu(t)}^{\mu(t)}
(y(t)))}{\partial x} \left( \Phi_{\nu(t)}^{\mu(t)} \right)_{\ast}
\dot y (t) = \frac{\partial \varphi_j (t,\Phi_{\nu(t)}^{\mu(t)}
(y(t)))}{\partial x}
 \left( \Phi_{\nu(t)}^{\mu(t)}
\right)_{\ast} a(y(t)),$$ $$ j{=}1,...,p,\; \; \; \; \; \; \;  \;
\; t\in[t_1{+}\sigma,T]$$ are equivalent to
$$ \frac{\partial {\tilde \psi}_j}{\partial y} (t,y(t)) \dot y (t)
= \frac{\partial {\tilde \psi}_j}{\partial y} (t,y(t)) a(y (t)),
\; \; \; \; j=1,...,p. \n{L17}$$ Finally, by the construction,
$\mu(T)=0,$ which yields $y(T)= \Phi_{\nu(T)}^{-\mu(T)}
(x_{\varkappa}(T)) =x_{\varkappa}(T)=x^T,$ and $\mu(t_1{+}\sigma)
= 0,$ which yields: $y(t_1{+}\sigma)=
\Phi_{\nu(t_1{+}\sigma)}^{-\mu(t_1{+}\sigma)}
(x_{\varkappa}(t_1{+}\sigma)) =x_{\varkappa}(t_1{+}\sigma).$

Combining this with (\ref{L4a}) and (\ref{L12_b}), we obtain:
$y(t_1{+}\sigma) \in \overline{\Omega_{\varepsilon_2}
(z(t_1{+}\sigma))},$ which yields condition 3) of lemma 5.1. Let
$(t,x) \mapsto \psi(t,x) = (\psi_1(t,x),...,\psi_n(t,x)) \in \rn$
be any map of class $C^1$  defined in some neighborhood ${\tilde
G}{\subset} [t_1{+}\sigma,T]{\times}\rn$ of the set $\{
(t,x)\in[t_1{+}\sigma,T] \times \rn \; | \; \; x=y(t),\; \;
t\in[t_1{+}\sigma,T] \}$ such that, for every fixed $t\in I,$ we
have $\psi(t,y(t)) = 0,$ and the map $x\mapsto
(\psi_1(t,x),...,\psi_n(t,x))$ defines canonical coordinates for
system (\ref{r3}) in the neighborhood ${\tilde G}_t:=\{ z \in \rn
\; | \; (t,z){\in} {\tilde G} \}$ of $y(t).$ (Since $y(\cdot)$ is
already defined, we can easily pick such a map following the same
pattern as that proposed in section 4 when proving theorem 3.1 as
a corollary of theorem 4.1. For this, we should consider the map
$(t_1,...,t_n) \mapsto (\Phi_{{\omega}_n}^{t_n} \circ
\Phi_{{\omega}_{n{-}1}}^{t_{n{-}1}}\circ... \circ
\Phi_{{\omega}_1}^{t_1}) (y(t))$ for each fixed $t{\in}
[t_1{+}\sigma,T],$ where ${\omega}_i(\cdot),$ $i{=}1,...,n$ are
any vector fields on $\rn$ such that $\Delta_i (x) = {\rm span }
\{{\omega}_{n-i} (x),{\omega}_{n-i+1} (x),..., {\omega}_n (x)\}$
for all $x {\in} \rn, $ $i{=}0,...,n{-}1.$ This map is a local
diffeomorphism in some neighborhood of $t_i{=}0,$ $i{=}1,...,n,$
and the inverse map $x\mapsto \psi (t,x)$ defines canonical
coordinates in a neighborhood of $y(t)$ for every fixed $t
\in[t_1{+}\sigma , T],$  and satisfies the condition $\psi(t,y(t))
= 0 \in \rn,\;$ $t \in[t_1{+}\sigma,T]).$ Then, since the
definition of the Lie derivative is coordinate-free, we obtain
from (\ref{L17})
$$    \frac{\partial\psi_j (t,y(t))}{\partial x} \dot y(t) =  \frac{\partial\psi_j (t,y(t))}{\partial x} a(y(t)), \;\; \; \; \;  j=1,\ldots,p, \; \; \; t{\in}[t_1{+}\sigma,T] $$
$$    \frac{d\psi_j (t,y(t))}{dt} =0, \; \; \; \; \; \; \; \; \; \; \; \; \; \;  j=1,\ldots,n, \; \;
\;t{\in}[t_1{+}\sigma,T] $$ which yields conditions 1), and 2) of
lemma 5.1.

Thus $y(\cdot),$ and $\psi(\cdot,\cdot) {=}
(\psi_1(\cdot,\cdot),...,\psi_n(\cdot,\cdot))$ satisfy conditions
1), 2), 3) of lemma 5.1. The proof of lemma 5.1 is complete.


\vskip10mm
\begin{center}
{\bf 7. Appendix}
\end{center}

\begin{center}
{\bf 7.1. Proof of lemma 6.1.}
\end{center}

From conditions (A),(B) it follows that for each $(t,x)$ in $\T=\{
(t,x)\;  |\; \; t{\in}I, \; \; |x-z(t)| \le \varepsilon_2\}$ (see
(\ref{L1add})) there exist a smooth vector field $\nu_{t,x}
(\cdot) \in \Delta_{n-p-1} (\cdot),$ and a point $y =
\Phi_{\nu_{t,x}}^{-\mu_{t,x}} (x)$ with some $\mu_{t,x}{\ge}0$
such that
$$\frac{\partial \varphi_p}{\partial t} (t,x ) + \frac{\partial
\varphi_p}{\partial x} (t,x )
\left(\Phi_{\nu_{t,x}}^{\mu_{t,x}}\right)_{\ast}
a(\Phi_{\nu_{t,x}}^{-\mu_{t,x}}(x)) = 0.$$ (Actually, there is a
finite sequence $\nu_{t,x}$ of vector fields from $\Delta_{n-p-1}
(\cdot)$ and the corresponding sequence, of nonnegative numbers,
which we denote by $\mu_{t,x},$ satisfying the above-mentioned
equality, as in the formulation of lemma 6.1. However, to make the
arguments clearer, we again assume (without loss of generality)
that each of these sequences consists on one element only, then we
denote these elements by $\nu_{t,x},$ and $\mu_{t,x}$
respectively. The proof for the general case is the same). Let
$\T_{t,x}$ be an open set of the form $\T_{t,x} =]t{-}\alpha,
t{+}\alpha [ \times \Omega_{\beta} (x)$ with some $\alpha{>}0,$
$\beta{>}0$ such that for its closure ${\overline{\T}}_{t,x} =
[t{-}\alpha,t{+}\alpha]\times\ \overline{\Omega_{\beta} (x)}$ we
get $\Phi_{\nu_{t,x}}^{-\mu} (\overline{\Omega_{\beta} (x)})
\subset \D_{\nu_{t,x}}$ for all $\mu {\in} [0,\mu_{t,x}],$ and
$$|\frac{\partial \varphi_p}{\partial t} (s,z ) + \frac{\partial
\varphi_p}{\partial x} (s,z )
\left(\Phi_{\nu_{t,x}}^{\mu_{t,x}}\right)_{\ast}
a(\Phi_{\nu_{t,x}}^{-\mu_{t,x}}(z))| < \delta \; \; \; \;  \; \;
\mbox{ for all } (s,z) \in {\overline \T}_{t,x}, \n{app27}$$ where
$\delta$ is defined in (\ref{rr24}). Since $\T {\subset}
\bigcup\limits_{(t,x){\in}\T} \T_{t,x},$ and $\T$ is a compact
set, there exists a finite open subcovering $\{ \T_r {=}
\T_{t_r,x_r} {\}}_{r=1}^{r_0}$ such that $\T \subset
\bigcup\limits_{r=1}^{r_0} {\T}_r.$ To simplify the notation, we
put by definition:
$$ \T_r : = \T_{t_r,x_r}, \; \; \nu_r (\cdot) : = \nu_{t_r,x_r}(\cdot), \; \;
\mu_r : = \mu_{t_r,x_r}, \; \; y_r:= \Phi_{\nu_r}^{-\mu_r} (x_r),
\; \; r=1,...,r_0. \n{app28}$$

Then $$|\frac{\partial \varphi_p}{\partial t} (t,x ) +
\frac{\partial \varphi_p}{\partial x} (t,x )
\left(\Phi_{\nu_r}^{\mu_r}\right)_{\ast}
a(\Phi_{\nu_r}^{-\mu_r}(x))| < \delta  \; \; \; \; \; \mbox{ for
all } (t,x) \in {\overline{\T}}_r \n{app29}$$

By definition, put:
$$L:=\frac{1}{2(M+1)}, \; \; \; \mbox{ where} \; \; M:= \max\{ \; \; |\digamma (t,x,v)|\; \; | \; \; \; v =  \frac{\partial \varphi_p (t,x )}{\partial t}  + $$
$$ +  \frac{\partial \varphi_p (t,x )}{\partial x}
(\Phi_{\nu_{r}}^{\mu})_{\ast} a(\Phi_{\nu_{r}}^{-\mu}(x)), \; \;
\; \; (t,x) \in {\overline{\T}}_r, \; \; 0 {\le} \mu {\le}
\mu_i,\; \; i{=}1,...,r_0 \} \n{app32}$$ Given any $(t,x) $ in
$\bigcup\limits_{r=1}^{r_0} {\T}_r,$ define $\theta_{t,x}(\cdot)$
and $\tau_{t,x} (\cdot)$ as follows
$$\theta_{t,x}(z):=t+\alpha_{t,x}^0 - L|z-x|; \; \; \; \tau_{t,x}(z):=t-\alpha_{t,x}^0 + L|z-x| \; \; \; \mbox{ for all } z\in \rn,$$
where $\alpha_{t,x}^0>0$ is a small positive number such that for
each set
$$ S_{t,x}:= \{ (s,z) \in I \times \rn \; | \; \; \tau_{t,x}(z) <s< \theta_{t,x}(z)\}$$
there is $r{=}r(t,x){\in}\{1,...,r_0\}$ such that $S_{t,x} \subset
\T_{r(t,x)}.$ Since $\T {\subset} \bigcup\limits_{(t,x){\in}\T}
S_{t,x},$ and $\T$ is compact, there is a finite open covering
$\{S_{t_m,x_m} {\}}_{m=1}^{m_0}$ of $\T:$ $\T \subset
\bigcup\limits_{m{=}1}^{m_0}S_{t_m,x_m}. $ Again, in order to
simplify the notation, we put
$$ \theta_m(\cdot)=\theta_{t_m,x_m}(\cdot); \; \; \tau_{m}:=\tau_{t_m,x_m} (\cdot); \; \;
S_m :=\{ (s,z) \in I \times \rn \; | \; \; \tau_m(z) \le s \le
\theta_m(z)\};$$
$$\T_m:=\T_{r(m)}; \; \; \nu_{m} (\cdot):=\nu_{r(t_m,x_m)} (\cdot);\; \; \mu_m:=\mu_{r(t_m,x_m)}\;
 \; z_m:=x_{r(t_m,x_m)} \n{app33} $$
Thus, from now on, we deal with notation (\ref{app33}), and
notation (\ref{app28}) is no longer valid. Let us remark that each
$\theta_m(\cdot)$ and each $\tau_m(\cdot)$ satisfy the global
Lipschitz condition $$\forall y{\in}\rn \; \; \; \forall z{\in}\rn
\; \; \; |\tau_j(y) -\tau_j (z)| \leq L\;  |y-z|; \; \; \;
|\theta_j(y)-\theta_j(z)| \leq L\;  |y-z|, \; \; \;
j{=}1,...,m_0.$$

 Let ${\Xi}$ be the system of all the sets given by
$$\Sigma_{\Theta(\cdot),\vartheta(\cdot),A_{\Theta},A_{\vartheta}}
:= \{ (s,z) \in \rr \times \rn \; | \; \; \vartheta(z) \le s \le
\Theta (z)\} \setminus \{ (s,z) \in \rr \times \rn \; | \; \; (s=
$$
$$= \vartheta(z), \; z\in A_{\vartheta}) \mbox{ or }   (s= \Theta(z),
\; z\in A_{\Theta}) \} \n{app34}
$$ where $\Theta (\cdot),$ and $\vartheta(\cdot)$ run through the
set of all the functions of class
 $C(\rn; I)$ such that, for all
$(y, z) {\in} \rn{\times}\rn,$
 $$ | \Theta(y) {-}\Theta(z)| {\leq} L |y{-}z| \; \; \;  \mbox{ and } \; \; \;   |\vartheta(y) {-} \vartheta(z) | {\leq} L  |y{-}z|,
 \;  \; \; \mbox{ for all } \; \;   y {\in} \rn, \;  z {\in} \rn,  \n{app35}$$
(with $L$ defined in (\ref{app32})) and $A_{\Theta} {\subset}
\rn,$ $A_{\vartheta} {\subset} \rn,$ run through the set of all
subsets of $\rn.$

Note that, if ${\vartheta}_j (\cdot),$ $j {=} 1,...,N,$ are some
functions of $\rn $  to $\rr$ such that
$$ \forall y \in \rn \;\;  \forall z \in
\rn \; \; \; \; \; \;\; \; \;    | {\vartheta}_j (y) -
{\vartheta}_j (z)| \leq
 L\;  |y-z|, \; \;j = 1,...,N,    $$
then, we obtain: $$ \forall y \in \rn \;\; \forall z \in \rn \;\;
\; \;  | \max\limits_{j=1,...,N} \{ {\vartheta}_j (y) \} -
\max\limits_{j=1,...,N} \{ {\vartheta}_j (z) \}| \leq
 L\;  |y-z|,   $$ and
$$ \forall y \in \rn \;\; \forall z \in \rn \;\;
\; \;  | \min\limits_{j=1,...,N} \{ {\vartheta}_j (y) \} -
\min\limits_{j=1,...,N} \{ {\vartheta}_j (z) \}| \leq
 L\;  |y-z|.   $$

 Therefore, it is easy to prove that $\Xi$ satisfies the following
 conditions: (a) $\emptyset {\in} {\Xi};$ (b) for each ${\Sigma}'
{\in} {\Xi},$ and each ${\Sigma}'' {\in} {\Xi},$
 we have ${\Sigma}'{\bigcap}{\Sigma}'' \in  {\Xi};$ and (c),
 for every $\Sigma {\in} {\Xi}, $ and every
$\Sigma_1 {\in} {\Xi},$ if $\Sigma_1 {\subset} \Sigma,$ then there
exists a finite sequence $\{ \Sigma_q {\}}_{q=2}^{q_0} \subset
{\Xi} $ of sets from ${\Xi}$ such that $\Sigma =
\bigcup\limits_{q=1}^{q_0} {\Sigma}_q,$ and ${\Sigma}_i {\bigcap}
{\Sigma}_j = \emptyset$ for all $i{\not=}j,$ $\{i,j\} {\subset}
\{1,...,q_0\}.$ (In other words, $\Xi$ is a "semiring" of sets  -
see \cite{kolmogorov}). Since every $S_m,$ $m{=}1,...,m_0$ is an
element of $\Xi,$ we get the existence of a finite sequence $\{
\Sigma_l {\}}_{l=1}^{l_0}$ of sets $\Sigma_l =
\Sigma_{\Theta_l(\cdot),\vartheta_l(\cdot),A_{\Theta_l},A_{\vartheta_l}}
\in \Xi$ (see (\ref{app34}), (\ref{app35})) such that first,
$\Sigma_{l'} \bigcap \Sigma_{l''} = \emptyset$ for all $l' \not=
l''$ in $\{1,...,l_0 \};$ second, $\T \subset
\bigcup\limits_{l{=}1}^{l_0}\Sigma_l
=\bigcup\limits_{m{=}1}^{m_0}S_m ;$ and, third, for each $l$ in
$\{1,...,l_0 \}$ there is $m(l) $ in $\{ 1,...,m_0\}$ such that
$\Sigma_l \subset S_{m(l)}.$

For each $l\in\{ 1,...,l_0\},$ we define the feedback control
$v_l(t,x)$ in $\T_{m(l)}$ as follows
$$ v_l(t,x) {=}  \frac{\partial \varphi_p (t,x )}{\partial t}  {+} \frac{\partial
\varphi_p (t,x )}{\partial x}
\left(\Phi_{\nu_{m(l)}}^{\mu_{m(l)}}\right)_{\ast}
a(\Phi_{\nu_{m(l)}}^{-\mu_{m(l)}}(x)) \; \;
 \; \; \mbox{ whenever }  (t,x) \in \T_{m(l)}, \n{app36_a}$$
and then we define the following (discontinuous!) feedback control
$v=v(t,x)$ in $\Sigma:=\bigcup\limits_{l=1}^{l_0} \Sigma_l$ for
system(\ref{rr21}) $$v(t,x):=v_l(t,x) \; \; \; \; \; \mbox{
whenever } (t,x) \in \Sigma_l, \; \; l=1,...,l_0. \n{app36}$$

Then, the following statement holds.

{\bf Lemma 7.1} {\it There are a unique trajectory $x(\cdot) \in
C(I;\rn)$ of system (\ref{rr21}) with the feedback law $v=v(t,x)$
given by (\ref{app36}), and with the initial condition $x(T)=x^T,$
a unique finite sequence of indices $\{ l_j {\}}_{j=1}^{N} \subset
\{ 1,...,l_0\}$ and a unique sequence
$T=\tau_{1}^{\ast}>\tau_2^{\ast}>...>\tau_{N}^{\ast}>\tau_{N{+}1}^{\ast}
=t_1$ such that  $l_i \not= l_j$ for all $i\not=j$ in $\{1,...,N
\}$ and

1) $\dot x (\cdot)$ is defined and continuous at each t in $I
\setminus  \{ \tau_2^{\ast},....,\tau_N^{\ast} \}$ and
$$(t,x(t)) \in \T \; \; \; \; \; \;  \mbox{ and }  \; \; \;   |v(t,x(t))| < \delta   \n{app37}  $$
 for all $t \in I.$

2) For every j=1,...,N we obtain $$ (t,x(t)) \in \Sigma_{l_j}; \;
\; \; \dot x (t) = \digamma (t,x(t),v(t,x(t))) \; \; \; \; \;
\mbox{ for all } \; \; \; t{\in}]\tau_{j{+}1}^{\ast},
\tau_j^{\ast}[ \n{app38}$$
$$ \tau_j^{\ast} = \Theta_{l_j} (z(\tau_j^{\ast})); \; \; \;\tau_{j+1}^{\ast} = \vartheta_{l_j} (z(\tau_{j+1}^{\ast})); \n{app39}$$
where $\digamma$ is defined in (\ref{rr31}), and
$\Theta_l(\cdot),$ $\vartheta_l (\cdot)$ are given in the
definition of $\Sigma_l.$ }

\vskip5mm

{\bf Proof of lemma 7.1.} We will prove the existence and the
uniqueness of $x(\cdot)$ and the corresponding $\{ \Sigma_{l_i}
{\}}_{i=1}^{N}$ by the induction over $i{=}1,...,N.$ In addition,
we will prove (by induction on $i{=}1,...,N$) that the trajectory
$x(\cdot)$ and the functions $(t,x){\mapsto} s_l(t,x),$ and
 $(t,x) {\mapsto }t_l (t,x),$ given by
 $$ s_l (t,x) =t - \vartheta_l(x), \; \; \; \; \; \; \; \;  t_l (t,x) =t - \Theta_l(x), \; \; \; \; \; \; \; \; (t,x)\in \rr\times\rn \n{app40}$$
 satisfy the conditions
 $$\frac{3(t-\tau)}{2} \geq s_l(t,x(t)) - s_l (\tau,x(\tau)) \geq \frac{t-\tau}{2}; \;  \; \; \; \;
  \frac{3(t-\tau)}{2} \geq t_l(t,x(t)) - t_l (\tau,x(\tau)) \geq $$ $$ \geq \frac{t-\tau}{2} \;\;  \; \; \;\; \; \;  \mbox{ for all } t>\tau, \; \;
   \{t,\tau\} \subset [\tau_i^{\ast},T], \; \; l=1,...,l_0. \n{app41}$$

For $i=1,$ we put by definition: $\tau_i^{\ast}=\tau_1^{\ast} =T.$
At this stage, we have empty set of $l_j,$ $\Sigma_{l_j},$ and
empty set of equalities (\ref{app38}), (\ref{app39}),
(\ref{app41}), $j{=}1,...,i{-}1,$ and the algorithm for getting
$l_1,$ $\Sigma_{l_1},$ and $\tau_2^{\ast}$ is the same as in the
general case $i {\ge} 1,$ which we consider now.

Assume that, for some $i \ge 1,$ there are a unique sequence
$\{l_j {\}}_{j=1}^{i-1}$ (for $i=1,$ the sequence is empty - see
the previous paragraph) such that $l_j \not= l_q$ for all $j \not=
q$ in $\{ 1,..., i{-}1 \},$ a unique sequence $T{=}\tau_{1}^{\ast}
{>} \tau_{2}^{\ast}{>}...{>}\tau_i^{\ast} {\ge} t_1,$ and a
trajectory $x(\cdot) {\in} C([\tau_i^{\ast},T]; \rn)$ such that
(\ref{app37}), and (\ref{app41})  hold for all $t\in
[\tau_{i}^{\ast}, T];$  and (\ref{app38}), (\ref{app39}) hold for
all $j{=}1,...,i{-}1;$  $\dot x(\cdot)$ being well-defined and
continuous on $[{\tau}_i^{\ast},T]\setminus
\{\tau_{2}^{\ast},...,\tau_{i}^{\ast} \}$ (again, for $i=1,$ we
deal with  empty set of equalities (\ref{app38}), (\ref{app39}),
(\ref{app41}), and (\ref{app37}) becomes trivial).

If $\tau_i^{\ast} {=} t_1,$ the proof is complete, therefore we
assume that $\tau_i^{\ast} {>} t_1.$ By the induction hypothesis,
$(\tau_i^{\ast}, x(\tau_i^{\ast})) \in \T;$ then, since $\T
\subset {\rm int} \left( \bigcup\limits_{l=1}^{l_0} \Sigma_l
\right)$ by the construction, we get the existence of $\alpha_0>0$
such that $(\tau_{i}^{\ast} {-} s , x(\tau_{i}^{\ast})) {\in}
\bigcup\limits_{l=1}^{l_0} \Sigma_l$ for all $s \in]0, \alpha_0].$
Since $\Sigma_{l'}\bigcap \Sigma_{l''} = \emptyset$  for all $l'
\not= l'',$ there are a unique $l_i {\in} \{1,...,l_0 \}$ and a
unique ${\overline \tau} \in [t_1,\tau_i^{\ast}[$ such that $
]{\overline\tau}, \tau_i^{\ast} [ \times \{ x(\tau_i^{\ast}) \}
\subset \Sigma_{l_i},$ and ${\overline\tau}= \vartheta_{l_i}
(x(\tau_i^{\ast})),$ $\tau_i^{\ast}= \Theta_{l_i}
(x(\tau_i^{\ast})).$ Since $v_{l_i} (\cdot,\cdot)$ is well-defined
in $\T_{m(l_i)}$ (see (\ref{app36_a})), and $(\tau_i^{\ast},
x(\tau_i^{\ast})) {\in} \Sigma_{l_i} \subset \T_{m(l_i)}, $  the
trajectory $t {\mapsto}
x(t,\tau_i^{\ast},x(\tau_i^{\ast}),v_{l_i}(\cdot,\cdot))$ is
well-defined on some maximal $]\overline{s}, \tau_i^{\ast}]$ such
that $x(t,\tau_i^{\ast},x(\tau_i^{\ast}),v_{l_i}(\cdot,\cdot)) \in
\T_{m(l_i)}$ for all $t \in ]\overline{s},\tau_i^{\ast}[. $ To
simplify the notation, we put $x_{l_i}(t):= x(t,\tau_i^{\ast},
x(\tau_i^{\ast}), v_{l_i} (\cdot,\cdot))$ for $t\in]\overline{s},
\tau_i^{\ast}].$ Let us remark that from
(\ref{rr31}),(\ref{app32}), (\ref{app35}) (and from the inclusions
$\Sigma_l {\in} \Xi,$ $l{=}1,...,l_0$) it follows that, for every
$l' {\in} \{1,...,l_0 \}$ such that $(\tau_i^{\ast},
x(\tau_i^{\ast})) {\in} \T_{m(l')}$ the trajectory $t\mapsto
x_{l'}(t):= x(t,\tau_i^{\ast}, x(\tau_i^{\ast}), v_{l'}
(\cdot,\cdot))$ satisfies the conditions
$$ \frac{3(t-\tau)}{2} \geq s_l(t,x_{l'}(t)) -  s_l(\tau,x_{l'}(\tau)) \geq \frac{(t-\tau)}{2}$$
$$ \frac{3(t-\tau)}{2} \geq t_l(t,x_{l'}(t)) -  t_l(\tau,x_{l'}(\tau)) \geq \frac{(t-\tau)}{2},
 \; \; \; t>\tau,\; \; l=1,...,l_0 \n{app42}$$
In particular, (\ref{app42}) holds for $l'{=}l_i$ and for all $t,$
$\tau$ in $]\overline{s}, \tau_i^{\ast}[,$ which implies that $t
\mapsto s_l (t, x_{l_i} (t)),$ and  $t \mapsto t_l (t, x_{l_i}
(t))$ are strictly increasing on $]\overline{s}, \tau_i^{\ast}[.$
From the induction hypothesis (see (\ref{app39}) for
$j=1,...,i{-}1$) it follows that
$t_l(\tau_i^{\ast},x(\tau_i^{\ast})) = 0,$
$s_l(\tau_i^{\ast},x(\tau_i^{\ast})) > 0.$ Since $s_l(t,x)>0,$ and
$t_l(t,x)<0$ are equivalent to $t>\vartheta_l(x)$ and
$t<\Theta_l(x)$ respectively, we obtain that
$(t,x(t,\tau_i^{\ast}, x(\tau_i^{\ast}), v_{l_i} (\cdot,\cdot)))$
belongs to $\Sigma_{l_i}$ for all $t{\in}]\tau_i^{\ast} {-}
\alpha_0, \tau_i^{\ast}[$ with some $\alpha_0{>}0,$ and moreover
there is a unique $\tau_{i{+}1}^{\ast} \in ]\overline{s},
\tau_i^{\ast}[$ such that $t_{l_i} (t,x_{l_i}(t)) < 0,$ $s_{l_i}
(t,x_{l_i} (t)) > 0$ for all  $t \in ]\tau_{i+1}^{\ast},
\tau_i^{\ast} [,$ and $t_{l_i} (t,x_{l_i}(t)) < 0,$ $s_{l_i}
(t,x_{l_i} (t)) < 0$ for all $t \in ]\overline{s},
\tau_{i+1}^{\ast} [,$ which implies that
$$(t,x_{l_i} (t)) \in \Sigma_{l_i}, \; \; \; \; \dot x_{l_i} (t) = \digamma (t,x_{l_i} (t), v_{l_i} (t,x_{l_i}(t)))
\; \; \; \;\; \; \; \; \; \; \;  \mbox{ for all } \; \;  t \in
]\tau_{i+1}^{\ast}, \tau_i^{\ast} [ \n{app43}$$
$$(t,x_{l_i} (t)) \notin \Sigma_{l_i} \; \; \; \;\; \; \; \; \; \; \; \;  \mbox{ for all } \; \;  t \in ]\overline{s}, \tau_{i+1}^{\ast}[  \n{app44}$$
Taking into account that $x_{l_i} (\tau_i^{\ast}) {=}
x(\tau_i^{\ast})$ by the construction, we obtain from
(\ref{app43}) and from the induction hypothesis that our $l_i,$
$\Sigma_{l_i},$ $\tau_{i{+}1}^{\ast},$ and $x(t)$ given by
$x(t):=x_{l_i} (t)$ for all $t\in [\tau_{i{+}1}^{\ast},
\tau_i^{\ast}]$ satisfy (\ref{app37})-(\ref{app39}), and
(\ref{app41}) for all $j=1,...,i{-}1,i.$

The uniqueness of $l_i,$ $\Sigma_{l_i},$ $\tau_{i{+}1}^{\ast},$
and $x(\cdot) $ on $[\tau_{i+1}^{\ast}, \tau_i^{\ast} ]$ such that
(\ref{app43}) holds with some ${\tau}_{i+1}^{\ast}$
(${\tau}_{i+1}^{\ast}{<} \tau_i^{\ast}$) follows from
(\ref{app42}), which is true for each $l'$ in $\{1,...,l_0 \}$
such that $(\tau_i^{\ast},x(\tau_i^{\ast})) \in \T_{m(l')}.$
Indeed, if there are  $l' {\not=} l_i$ in $\{1,...,l_0 \},$
$\overline{\tau}_{i+1}^{\ast}{<} \tau_{i}^{\ast} ,$ and
${\overline x} (\cdot) \in C([\overline{\tau}_{i+1}^{\ast},
\tau_{i}^{\ast}]; \rn)\bigcap C^1(]\overline{\tau}_{i+1}^{\ast},
\tau_{i}^{\ast}[; \rn) $ such that $${\overline x}
(\tau_{i}^{\ast}) {=}{x} (\tau_{i}^{\ast}); \; \; \; \; \; \; \;
{\overline x}(t) \in \Sigma_{l'} \; \; \; \; \; \; \mbox{ for all
} \; \;   t \in ]\overline{\tau}_{i+1}^{\ast}, \tau_{i}^{\ast}[$$
$$ \dot{\overline x} (t) = \digamma (t,{\overline x} (t), v(t,{\overline x} (t) )) =
 \digamma (t,{\overline x} (t), v_{l'}(t,{\overline x} (t) )), \; \; \;
  t\in ]\overline{\tau}_{i+1}^{\ast}, \tau_{i}^{\ast}[$$
$$ \vartheta_{l'} (\overline{\tau}_{i+1}^{\ast}, {\overline x}(\overline{\tau}_{i+1}^{\ast})) = \overline{\tau}_{i+1}^{\ast},  \n{app45} $$
then we get from (\ref{app42}) and from the definition of $l_i:$
$s_{l_i}(t,{\overline x}(t))>0,$ $t_{l_i}(t,{\overline x}(t))<0$
for $t \in ]\overline{\tau}, \tau_{i}^{\ast}[$ with some
$\overline{\tau} \in ]\overline{\tau}_{i+1}^{\ast},
\tau_{i}^{\ast}[,$ which yields: $(t,\overline{x}(t)) \in
\Sigma_{l'} \bigcap \Sigma_{l''}$ for all $t \in ]\overline{\tau},
\tau_{i}^{\ast}[.$ Since $\Sigma_{l'} \bigcap \Sigma_{l_i} =
\emptyset$ for each $l' \not= l_i,$ we obtain that $l' = l_i,$
which proves the uniqueness of $l_i$ and of $\Sigma_{l_i}.$ This,
in turn, implies that $\overline{x} (t) = x(t)$ for all $t \in
[\hat{\tau}, \tau_i^{\ast}],$ where $\hat{\tau} := \max \{
{\overline{\tau}}_{i{+}1}^{\ast}, \tau_{i{+}1}^{\ast} \}.$ The
function $s_{l_i}(t) = t -\vartheta_{l_i}(t,x(t))$ is strictly
increasing on $[{\overline{\tau}}_{i{+}1}^{\ast}, \tau_i^{\ast}]
\bigcup [{{\tau}}_{i{+}1}^{\ast}, \tau_i^{\ast}];$ therefore, from
the equalities
$\vartheta_{l_i}({\overline{\tau}}_{i{+}1}^{\ast},{\overline
x}({\overline{\tau}}_{i{+}1}^{\ast})) {=}
{\overline{\tau}}_{i{+}1}^{\ast}, $
$\vartheta_{l_i}({\tau}_{i{+}1}^{\ast},
x({\tau}_{i{+}1}^{\ast})){=} {\tau}_{i{+}1}^{\ast},$ it follows
that $\hat{\tau} {=} \overline{\tau}_{i{+}1}^{\ast} {=}
{\tau}_{i{+}1}^{\ast},$ which proves the uniqueness of
${\tau}_{i{+}1}^{\ast},$ and the uniqueness of $x(\cdot)$ on
$[{\tau}_{i{+}1}^{\ast}, {\tau}_{i}^{\ast}].$ Let us remark that
$(t,x(t)) \in \T$ for all $t \in [{\tau}_{i{+}1}^{\ast},
{\tau}_{i}^{\ast}],$ which follows from the inclusion $(t,x(t))
\in \T$ for $t \in [{\tau}_{i}^{\ast},T]$ (see the induction
assumption) and from the choice of $\delta$ (see
(\ref{rr24}),(\ref{app27})). Finally, (\ref{app41}), and
(\ref{app42}) with $l'=l$ yield
$$ \frac{3(t-\tau)}{2}\geq s_l(t,x(t))-s_l(\tau,x(\tau)) \geq \frac{(t-\tau)}{2}; \; \; \;  \;
 \frac{3(t-\tau)}{2}\geq t_l(t,x(t))-t_l(\tau,x(\tau)) \geq $$
 $$\geq  \frac{(t-\tau)}{2}\; \; \; \; \; \; \; \; \mbox{ for all } t>\tau, \; \;
 \; \{ t, \tau\} \subset [\tau_{i{+}1}^{\ast}, T], \; \; \;
 l=1,...,l_0, $$ which implies that $s_{l_i}(t,x(t)) {>} 0
 {>}s_{l_j}(t,x(t)),$  $j{=}1,...,i{-}1,$ and therefore $l_i \not= l_j$
 for all $j{=}1,...,i{-}1;$ hence $l_j \not=l_q$ for all $q\not=j,\; $
 $\{q,j \} \subset \{1,...,i\}.$ The construction of $\Sigma_{l_j},$
 $\tau_{j+1}^{\ast},$  and $x(t),$ $t\in [\tau_{j+1}^{\ast},
 \tau_{j}^{\ast}]$ is complete.

  Thus, we obtain by the induction
 sequences $\{ l_j \},$  $\{ \Sigma_j \},$ $\tau_j^{\ast},$ and
 the trajectory $x(\cdot)$ of system (\ref{rr21}) with the
 feedback control $v=v(t,x)$ given by (\ref{app36}) such that
 $x(T) = x^T,$ $\tau_{j{+}1}^{\ast} < \tau_{j}^{\ast};\; \; $
 $l_i \not= l_j$ for all $i \not= j,$ and
 (\ref{app37})-(\ref{app39}) hold for all $j=1,2,\ldots \; $ .  If $\tau_{i+1}^{\ast} = t_1$
 for some $i \in \nn,$ we put $N:=i,$ and the proof of lemma 6.1
 is complete. Otherwise, sequences $\{ \tau_i^{\ast} \},$  and $\{ l_i \}$
 are infinite, which is impossible because $\{ l_i \} \subset \{ 1,..., l_0
 \}$ and $l_i \not= l_j,$ whenever $i \not= j.$ The proof of lemma
 7.1. is complete.
Finally, we put $\T_i:=\T_{m(l_i)},$ and $v_i(t,x):=v_{l_i}
(t,x),$ for all $(t,x) \in \T_i:=\T_{m(l_i)},$ $i=1,...,N.$ Then
we obtain from lemma 7.1 that our sequences of $\T_i$ and
$v_i(\cdot,\cdot)$ satisfy all the conditions of lemma 6.1. The
proof of lemma 6.1. is complete.

\vskip10mm




\end{large}
\end{document}